\documentclass{article}

\usepackage{arxiv}
\usepackage{bold-extra}
\usepackage[english]{babel}
\usepackage{caption}
\usepackage{lscape} 
\usepackage{amsmath}
\usepackage{geometry}
\usepackage{amsthm}
\usepackage{adjustbox}
\usepackage{amsfonts}
\usepackage{mathtools}
\usepackage{verbatim}
\usepackage{algpseudocode,algorithm,algorithmicx}
\usepackage{mathtools}
\usepackage{comment}
\usepackage{verbatim, booktabs, dcolumn, setspace, fullpage}

\newtheorem{theorem}{Theorem}

\newtheorem{cor}{Corolary}

\newtheorem{lemma}{Lemma}

\newtheorem{remark}{Remark}
\usepackage{graphicx}
\usepackage{wrapfig}
\usepackage{caption}
\usepackage{mathtools}
\usepackage{float}

\allowdisplaybreaks

\usepackage{tikz}
\usetikzlibrary{matrix, positioning, arrows.meta, decorations.pathreplacing}

\usepackage[english,noautotitles-r]{SASnRdisplay} 
\usepackage[round]{natbib}
\usepackage{url}
\usepackage{subfigure}
\lstdefinestyle{r-output}{
style = r-style,
style = r-output-user,
}

\usepackage{url}            
\usepackage{booktabs}       
\usepackage{amsfonts}       
\usepackage{nicefrac}       
\usepackage{microtype}      
\usepackage{graphicx}
\usepackage{natbib}
\usepackage{doi}

\newcommand{\E}{\mathbf{E}}
\newcommand{\D}{\mathbf{Var}}
\newcommand{\Probb}{\mathbf{P}}

\newcommand{\Cov}{\mathbf{Cov}}

\title{Testing MCAR via covariances: Extending the U-statistic framework to partially observed variables}

\author{ \href{https://orcid.org/0000-0002-0460-400X}{\includegraphics[scale=0.4]{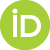}\hspace{1mm} Danijel G. Aleksi\' c} \\
	University of Belgrade\\
	Faculty of Organizational Sciences\\
    Department of Mathematics \\
	Belgrade, 11000, Serbia \\
	\texttt{danijel.aleksic@fon.bg.ac.rs} \\
}	

\date{}

\begin{document}	
\maketitle

\begin{abstract}
This paper presents a generalized version of a \textit{U}-statistics-based test for MCAR developed by \cite{aleksic2024novel}. The proposed test, similar to the original, evaluates the MCAR assumption by calculating and combining the covariances between response indicators and data variables. However, unlike the preceding version, it is capable of utilizing partially observed variables, resulting in a significantly larger class of detectable alternatives. Numerical results indicate that the improved test is well-calibrated, notably outperforming the well-known MCAR test developed by \cite{Little1988} used as a benchmark. For alternatives detectable by the original method, the improved test maintains comparable, although slightly lower, power, while consistently outperforming Little's test across all studied scenarios. For alternatives that were previously undetectable or marginally detectable, the novel test demonstrates the superior performance among the three methods. While the novel test shares the assumption of finite fourth moments of the data with Little's test, the results suggest it is more robust to this requirement, although both tests exhibit similar limitations.
    
    \vspace{6pt}\small
\end{abstract}

\section{Introduction}\label{sec:introduction}


Handling of missing data is a very important step in any data analysis. When data are missing completely at random (MCAR) \citep{LittleRubin1987}, and either the missingness rate is modest or the sample size is large, removing the incomplete observations from the data, so-called complete-case analysis, will suffice in the vast majority of cases. Moreover, there are many adaptations of existing tests that provide better type I error control and better power performance under MCAR data \citep{aleksic2023etAl, aleksic2024impute}. However, if data are not MCAR, complete-case analysis tends to be disastrous, producing biased estimates, shifted confidence intervals and tests with little to no type I error control. Knowing this, it is of essential importance to have a high-quality test for testing the MCAR assumption; consequently, there has been quite a lot of interest in testing whether it holds, or not.

\subsection{Paper aim and structure}

Recently, a \textit{U}-statistics-based test for MCAR developed by \cite{aleksic2024novel}, which performed better than the well-known MCAR test by \cite{Little1988}. However, as we will soon see, the test had a very narrow class of detectable alternatives. Overcoming that fundamental flow of that test is the main aim of this paper.

 The remaining of the paper is structured as follows. Subsection \ref{subsec:aleksic_test} presents, for readers' convenience, a brief overview of the original test. Section \ref{sec:improved_test} presents the improved version of the test, which is able to utilize most of the partially observed variables. The theoretical properties are also derived and some theoretical limitations are noted. The simulation study is given in Section \ref{sec:simulation}, where the improved test is compared to the old one in all of the originally studied settings, in all of which the old one performed better than Little's test. Some additional comparisons to Little's test are also given in settings where the old test could not perform at all due to its limitations.  Finally, some concluding remarks are given in Section \ref{sec:conclusion}.

\subsection{Tests of MCAR}

First results on MCAR testing were developed in the 1980's: for categorical data by \cite{fuchs1982maximum}, and for the Gaussian data by  \cite{Little1988}. As far as we know, Little's MCAR test, that is based on comparing the MLEs across missingness patterns, remains the most widely used test in practice. However, as we will see {in} Section \ref{sec:simulation}, it is highly sensitive to the assumption of underlying multivariate normality.

\cite{diggle1989testing} considered missing data in the context of {repeated measurements}, i.e. when a time-ordered sequence of measurements is made on some participants in {an} experiment. He considered a special case of missingness called a \textit{dropout}, where a sequence of measurements is terminated prematurely, and developed the class of procedures that test whether dropouts in the data occur randomly, in the sense that they are not related to any of the past measurements. {The methodology involves selecting a score function, where large values indicate rejection of the null hypothesis, and applying the normal approximation when feasible.} \cite{ridout1991testing} presented some improvements in terms of flexibility, utilizing logistic regression. 

{\cite{park1993atest} extended the Little's test to incomplete repeated categorical data.} Similarly, \cite{park1993btest} relied on Little's test to make a MCAR test for {repeated measurements}. Following the idea of Park and Davis, \cite{park1997test} constructed a MCAR test for the incomplete longitudinal data in the framework of generalized estimating equations. 

\cite{listing1998tests} developed a test for random dropouts in clinical trials by comparing the means of the individuals that stay, and those that drop out. 

Another test for the framework of generalized estimating equations, but for independent observations, came from \cite{chen1999test}, and {generalized} the original idea of \cite{Little1988}. \cite{qu2002testing} proposed a more unified generalized score-type test for ignorable missingness in longitudinal data.  

\cite{kim2002tests} studied tests based on weighted generalized least squares methods, and compared them to the likelihood-based tests, such as Little's test, in terms of type I error and power behavior in small sample sizes. The comparison {examines the} homogeneity of means and covariance matrices across missing data patterns.  

{A further} logistic regression-based testing procedure for MCAR tailored for medical longitudinal data was developed by \cite{fairclough2003design}. The main idea of {the} procedure is to study the dependence between the response indicators and the quality of life scores.

Although testing MCAR vs. MAR is not possible in the general case, since the data needed for that testing are missing, \cite{potthoff2006can} proposed {the test for \textit{MAR+} assumption, which is a special, testable case of MAR.}

The idea of testing MCAR by comparing the covariance matrices across the missing data patterns came with \cite{jamshidian2007testing}. \cite{jamshidian2008postmodeling} constructed a test for distinguishing MCAR from MNAR by noting that the maximum-likelihood estimates across random data subsamples will have the same asymptotic distribution under MCAR, but not under MNAR. 

\cite{fielding2009investigating} presented a real-data empirical comparison of four tests in the context of quality of life outcomes. {Specifically, the tests of} \cite{Little1988}, \cite{listing1998tests}, \cite{ridout1991testing}, and \cite{fairclough2003design} {were applied to several datasets to assess differences in inferential results.}

\cite{jamshidian2010tests} have considered {a} test for MCAR that relies on imputing the dataset and then conducting the complete-data procedures. The data are grouped by missingness patterns, and the variances across groups of data are then compared. \cite{jamshidian2013data} improved the results of \cite{jamshidian2008postmodeling} by approximating the asymptotic distribution {rather than using} the bootstrap method. 

\cite{lin2013probability} developed a probability based framework for testing MCAR that appeared comparable to Little's MCAR test in terms of power for a large number of studied scenarios. 

\cite{jamshidian2014examining} gave an overview of then {available} MCAR test that are based on either homogeneity of parameters or homogeneity of distributions across missingness patterns. Additionally, they proposed a novel nonparametric test for MCAR that is based on pairwise comparison of marginal distributions of the data, considering one variable at a time fixed. \cite{li2015nonparametric} also considered a nonparametric test for MCAR. The procedure first splits all of the data into categories by missingness patterns, and then uses Rizo-Sz\'ekely dissimilarity measure to compare distributions across patterns. The bootstrap algorithm is utilized afterwards {to approximate the $p$-value of the test.}

\cite{yuan2018missing} showed that, under normality, MLEs for different missingness patterns can converge to the same values, possibly not the true ones, even under MAR or MNAR. As a result, tests for MCAR based on comparing means and covariances across patterns cannot be safely used.

\cite{zhang2019unified} noted that most MCAR tests do not offer a method for a subsequent estimation once MCAR is rejected, and they presented a unified likelihood approach for both MCAR testing and subsequent estimation that appeared {to} behave well in the observed (although limited) scenarios. 

\cite{bojinov2020diagnosing} considered testing MAAR (Missing Always at Random), where response {mechanism} {does not depend on the data not only for the observed, but for any possible missingness pattern}. They noted that under certain regularity conditions, MAAR can be tested from the observed data only, and {proposed} three diagnostic procedures that rely on testing the dependence between response indicators and fully observed variables. 

\cite{spohn2021pklm} introduced the test that measures distributional differences across missing data patterns using Kullback-Leibler divergence. 

\cite{rouzinov2022regression} tested MCAR by fitting the linear regression model on the complete cases, and then comparing distributional differences of predicted values for missing and observed data. 

For the case of hidden Markov models, \cite{chassan2023test} developed a MCAR which does not require grouping the data by patterns, but are based on the estimates of conditional (given the latent state of the Markov chain) probabilities of missingness.

Lately, the measure of \textit{compatibility} was utilized by \cite{berrett2023optimal} for testing MCAR. Their key point is that there can be no test that can reject MCAR if the class of marginal distributions is compatible, i.e. they were successful in describing the exact class of non-detectable alternatives to MCAR. {They related the concept of compatibility testing to MCAR testing in the discrete case.} \cite{bordino2025tests} compared compatibility of covariance matrices across missing data patterns to construct a MCAR test for the incomplete data that do not need to be discrete. The formal definition of compatibility can be found in either of these works, and can be traced back to \cite{sklar1959fonctions} and the theory of copulas.

Dealing with functional data, we refer the reader to the recent test by \cite{ofner2025testing}.

Most of the existing statistical tests for MCAR are based on comparing some measure across different missing data patterns. To the best of our knowledge, there were no tests constructed using the rationale of checking the linear dependence between the response indicators and fully observed data columns. {This changed with the test by \cite{aleksic2024novel}, who utilized the non-degenerate U-statistics theory to combine the covariances between the complete data variables and the response idicators of the incomplete ones. As we will see from this paper, a significant gap for the improvement of that test was left open.}

\subsection{Original test}\label{subsec:aleksic_test}

In a recently developed test, \cite{aleksic2024novel} studied the random sample of $n$ independent copies of a random vector $\left( X^{(1)}, \dots, X^{(p)}, Y^{(1)}, \dots, Y^{(q)} \right)$, where he assumed that variables $X^{(1)},\dots,  X^{(p)}$ are completely observed, and variables $Y^{(1)}, \dots, Y^{(q)}$ are susceptible to missingness. More precisely, the sample can be written as
\begin{equation}\label{matrix_sample}
    \begin{bmatrix}
X^{(1)}_1 & \cdots & X^{(p)}_1 & Y_1^{(1)} & \cdots & Y_1^{(q)}  \\
X^{(1)}_2 & \cdots & X^{(p)}_2 & Y_2^{(1)} & \cdots & Y_2^{(q)} \\
\vdots & \ddots &\vdots & \vdots & \ddots & \vdots   \\
X^{(1)}_n & \cdots & X^{(p)}_n & Y_n^{(1)} & \cdots  & Y_n^{(q)}  
\end{bmatrix}.
\end{equation}
If, for $1 \leq i \leq n$ and $1 \leq j \leq q$ we introduce the response indicator $R^{(j)}_i$, which is equal to 1 if $Y^{(j)}_i$ is observed, and 0 if missing, the sample can be written in an expanded form:
\begin{equation}\label{expanded_sample_last}
    \begin{bmatrix}
X^{(1)}_1 & \cdots & X^{(p)}_1 & Y_1^{(1)} & \cdots & Y_1^{(q)} & R_1^{(1)} & \cdots & R_1^{(q)} \\
X^{(1)}_2 & \cdots & X^{(p)}_2 & Y_2^{(1)} & \cdots & Y_2^{(q)} & R_2^{(1)} & \cdots & R_2^{(q)} \\
\vdots & \ddots &\vdots & \vdots & \ddots & \vdots & \vdots & \vdots \\
X^{(1)}_n & \cdots & X^{(p)}_n & Y_n^{(1)} & \cdots  & Y_n^{(q)}  & R_n^{(1)} & \cdots & R_n^{(q)}
\end{bmatrix}.
\end{equation}

Aleksi\'c developed the MCAR test that was based on the fact that, when the equality
\begin{equation*}
    \Cov \left( X^{(u)}, R^{(v)} \right) = 0
\end{equation*}
does not hold for some $1 \leq u \leq p$ and $1 \leq v \leq q$, the MCAR assumption can be rejected. Having this, he proposed the test based on the unbiased estimates of $ \Cov \left( X^{(u)}, R^{(v)} \right)$:
\begin{equation}\label{Tn_u_v}
    T_{n,X}^{(u, v)} =  \frac{1}{n} \sum_{i=1}^n X_i^{(u)}R_i^{(v)} - \frac{1}{n(n-1)} {\sum_{i=1}^n \sum_{\substack{j = 1 \\ j \neq i}}^n} X_i^{(u)}R_j^{(v)} .
\end{equation}
If any of $T_{n,X}^{(u, v)}$ significantly deviates from zero, it is a strong sign that MCAR should be rejected. The main challenge was how to combine all of the $T_{n,X}^{(u, v)}$ into one test statistic. However, Aleksi\'c noted that $T_{n,X}^{(u, v)}$ is a difference of two non-degenerate $U$-statistics, and, utilizing their known properties, was able to prove that
\begin{equation}\label{big_convergence}
\left(  \sqrt{n} T_{n,X}^{(1,1)}, \dots, \sqrt{n} T_{n,X}^{(1, q)}, \sqrt{n} T_{n,X}^{(2,1)}, \dots, \sqrt{n} T_{n,X}^{(2,q)}, \dots, \sqrt{n}T_{n,X}^{(p,1)}, \dots, \sqrt{n}T_{n,X}^{(p,q)}  \right) \overset{D}{\to} \mathcal{N} (\boldsymbol{0}, \Sigma),
\end{equation}
as $n \to \infty$, where
\begin{equation}\label{sigma_compact}
    \Sigma = \Cov \left( \left( X^{(1)}, X^{(2)}\right)   \right) \otimes \Cov \left( \left( R^{(1)}, R^{(2)}\right)   \right),
\end{equation}
where $\otimes$ denotes the Kronecker product of (covariance) matrices.

Scaling by the standard estimate $\hat{\Sigma}^{-1}$ of  the inverse of $\Sigma$, under the assumption of finite fourth moments of the $X^{(1)}, \dots, X^{(p)}$, it was proven that, under the null hypothesis of MCAR,
\begin{equation}\label{An}
    A_n = n \left( T_{n,X}^{(1, 1)}, \dots, T_{n,X}^{(p, q)} \right) \boldsymbol{\hat{\Sigma}}^{-1} \left( T_{n,X}^{(1, 1)},  \dots, T_{n,X}^{(p, q)} \right)^T \overset{D}{\to} \chi^2_{pq},
\end{equation}
which can be used to create the rejection region for the test statistic $A_n$ of a novel test for MCAR.

The test appeared to perform much better than the well-known Little's MCAR test \citep{Little1988}, that is the most commonly used MCAR test in practice nowadays. The novel test performed significantly better in all of the studied cases, in both empirical size and power terms. 

Furthermore, for the special case of univariate nonresponse, it was shown that $A_n^2 = d^2$, where $d^2$ is the test statistic of Little's test.

\textit{However, the test left much to be desired}. The first major drawback of the test is that it can be used only on the dataset that has at least one complete column. The second and more important drawback is that it does not use the partially observed variables $Y^{(1)}, \dots, Y^{(q)}$ at all. That leads not only to the loss of power, but to potential inability to detect alternatives to MCAR where response indicators depend on $Y^{(1)}, \dots, Y^{(q)}$, but not on $X^{(1)}, \dots, X^{(p)}$. That setting is not very uncommon, so it is of a essential importance to address the issue. 

\begin{remark}
    The notation used throughout this subsection, and the rest of the paper, differs slightly from that of the original paper \citep{aleksic2024novel}. Additionally, the original approach rejected the MCAR hypothesis when the ``product theorem'' did not hold, i.e., when $\E(X)\E(R) \neq \E(XR)$. The test was constructed around the parameter $\E(X)\E(R) - \E(XR)$, which is equal to $-\Cov(X,R)$. This change in formulation has no impact on the performance of the test, as the covariance estimates are squared within the test statistic. The notation used here follows that established in \cite{aleksic2026PhD}.
\end{remark}


\section{Generalized test}\label{sec:improved_test}

Our main goal is to find a way to compute the statistic $T_n^{(u,v)}$ from \eqref{Tn_u_v} for each pair $(Y^{(u)}, R^{(v)})$, with $u \neq v$, in order to obtain an estimate of the covariance between them.
The subsequent goal is to use those statistics to extend the vector from \eqref{An} to include them, and, as a consequence, to expand the set of detectable alternatives of the test. Under the null hypothesis of MCAR data, it is reasonable to calculate it on those cases where $Y^{(u)}$ is observed. In that case, the statistic can be written as
\begin{align}
    T_{n, Y}^{(u, v)} &=  \frac{1}{\hat{n}^{(u)}} \sum_{i=1}^n Y_i^{(u)}R_i^{(u)}R_i^{(v)} - \frac{1}{\hat{n}^{(u)} \left(n - 1 \right)} {\sum_{i=1}^n \sum_{\substack{j = 1 \\ j \neq i}}^n} Y_i^{(u)} R_i^{(u)}R_j^{(v)}, \quad 
1 \leq u, v \leq q, \quad u \neq v, \label{Tn_u_v_trunc}
\end{align}
where $\hat{n}^{(u)} = \sum_{i=1}^n R_i^{(u)}$, $1 \leq u \leq q$, is the number of observed values in $Y^{(u)}$.

\begin{remark}
    We note that the form \eqref{Tn_u_v_trunc} is {a} generalization of \eqref{Tn_u_v}, since, for complete variables, all of the response indicators are equal to 1. 
\end{remark}

It is intuitive (and true) that, under MCAR data, a non-degenerate \textit{U}-statistic computed from complete cases, appropriately normalized, has the same asymptotic distribution as the one based on the complete sample. However, rigorous proof was anything but trivial \citep{aleksic2023etAl}. The complexity of the result is due to the fact that $\hat{n}^{(u)}$ is not constant, but random. Formally speaking, $T_{n, Y}^{(u, v)}$ is not a \textit{U}-statistic, but is asymptotically equivalent to one. In our case, which involves the difference between two test statistics and the joint distribution of such statistics, we should not expect the situation to be any simpler.

The main reason for introducing the statistic $ T_{n, Y}^{(u, v)}$ from \eqref{Tn_u_v_trunc} is that it serves as an approximation of an unbiased estimator of $\Cov \left( Y^{(u)}, R^{(v)}  \right)$, which is a measure of dependence between $Y^{(u)}$ and $R^{(v)}$. However, the estimate of any value proportional to it would also suffice. So, naturally, one could think of
\begin{equation}\label{Tn_u_v_hat}
    \hat{T}_{n, Y}^{(u, v)} =  \frac{1}{n} \sum_{i=1}^n Y_i^{(u)}R_i^{(u)}R_i^{(v)}- \frac{1}{n \left( n - 1 \right)} {\sum_{i=1}^n \sum_{\substack{j = 1 \\ j \neq i}}^n} Y_i^{(u)}R_i^{(u)}R_j^{(v)}. 
\end{equation}
as a modification of the test statistic. This choice seems appropriate, since we use deterministic $n$ instead of random $\hat{n}^{(u)}$, and 
\begin{align*}
    \hat{T}_{n, Y}^{(u, v)} = \frac{\hat{n}^{(u)}}{n}T_{n, Y}^{(u, v)},
\end{align*}
so it seems like Slutsky's theorem could be used as a final step of deriving the asymptotic behavior. However, a more fundamental problem lies beneath the surface. Under mutual uncorrelatedness of response indicators, it holds that
\begin{equation*}
    \E \left( \hat{T}_{n, Y}^{(u, v)}  \right) = \E\left( T_{n, Y}^{(u, v)} \right) \E\left( R^{(u)} \right) = \Cov \left( Y^{(u)}, R^{(v)} \right) \E\left( R^{(u)} \right),
\end{equation*}
which is equal to zero under MCAR. This makes $\hat{T}_{n,Y}^{(u,v)}$ appear to be a suitable choice, as its expected value is proportional to the target covariance. However, if response indicators exhibit any form of dependence among themselves, the statistic $T_{n,Y}^{(u,v)}$ from \eqref{Tn_u_v_trunc} is no longer a complete-case estimate under MCAR, since the data on which it is computed are not MCAR in that case. Indeed, the data here consist of realizations of the pair $(Y^{(u)}, R^{(v)})$, while complete cases are selected with respect to $R^{(u)}$. These response indicators corresponding to different variables need not be mutually independent for the data to be MCAR; the only assumption required by the test will be the \textit{absence of multicollinearity} among them, that is, their covariance matrix must be regular. 

If we assume MCAR and make no additional assumptions about the response indicators, it is straightforward to see that
\begin{align*}
    \E \left(\hat{T}_{n, Y}^{(u, v)}  \right) = \E \left( Y^{(u)} \right) \Cov \left( R^{(u)}, R^{(v)} \right),
\end{align*}
which does not necessarily equal zero under the null hypothesis, as it should be. In fact, $\hat{T}_{n, Y}^{(u, v)}$ is a \textit{U}-statistic that estimates $\Cov \left( Y^{(u)} R^{(u)}, R^{(v)} \right)$, which may or may not be close to desired $\Cov \left( Y^{(u)}, R^{(v)} \right) \, \E \left( R^{(u)} \right)$, depending on the internal covariance structure of the response indicators. In other words, the statistic $\hat{T}_{n, Y}^{(u, v)}$ can be interpreted as an indirect measure of the covariance between the incomplete variable and the response indicator, although with variable reliability. 

To be able to examine this issue in more detail, we first need asymptotic results for these statistics, which we state in the form of the following Lemma. The proof can be found in the Subsection \ref{lemma_YY_proof} of Appendix \ref{Appendix_proofs}.

\begin{lemma}\label{lemma_YY}
    Under MCAR, it holds that
    \begin{align}
    \lim_{n \to \infty}& \Cov \left( \sqrt{n} \hat{T}_{n, Y}^{(u, v)}, \ \sqrt{n} \hat{T}_{n, Y}^{(r, s)}  \right) = \E \left( Y^{(u)}Y^{(r)} \right) \cdot A + \E \left( Y^{(u)} \right)\E \left(Y^{(r)} \right) \cdot B ,\label{cov_YY}
\end{align}
where
\begin{align*}
    A &= \E \left(R^{(u)}R^{(v)}R^{(r)}R^{(s)}\right)  - \E \left( R^{(u)}R^{(v)}R^{(r)} \right) \E \left( R^{(s)}\right)\\
    & \hspace{60pt}-\E \left( R^{(u)}R^{(r)}R^{(s)} \right)\E \left( R^{(v)} \right)  + \E \left( R^{(u)}R^{(r)} \right) \E \left(R^{(v)} \right) \E \left( R^{(s)}\right)
\end{align*}
and
\begin{align*}
    B &= \E \left( R^{(u)}R^{(s)}  \right) \E \left( R^{(v)}\right)\E \left(R^{(r)} \right)  +   \E \left(R^{(v)}R^{(r)} \right) \E \left(R^{(u)} \right) \E \left( R^{(s)}\right)  + \E \left(R^{(v)}R^{(s)} \right) \E \left(R^{(u)} \right) \E \left( R^{(r)}\right) \\
    & \hspace{15pt}+2 \E \left( R^{(r)}R^{(s)} \right)  \E \left( R^{(u)}\right) \E \left( R^{(v)}\right) + 2 \E \left(R^{(u)}R^{(v)} \right)\E \left( R^{(r)}\right) \E \left(R^{(s)} \right) - \E \left( R^{(u)}R^{(v)}R^{(s)}\right) \E \left( R^{(r)} \right)\\
    &\hspace{15pt} - 4\E \left(R^{(u)} \right)\E \left(R^{(v)} \right)\E \left( R^{(r)} \right) \E \left(R^{(s)} \right) - \E \left( R^{(u)}R^{(v)} \right) \E \left( R^{(r)}R^{(s)} \right) - \E \left( R^{(v)}R^{(r)}R^{(s)}\right) \E \left(R^{(u)} \right).
\end{align*}
In particular, if either the incomplete variables have zero means or the response indicators are uncorrelated, it holds that
\begin{align}\label{cov_YY_omitted}
    \lim_{n \to \infty}& \Cov \left( \sqrt{n} \hat{T}_{n, Y}^{(u, v)}, \ \sqrt{n} \hat{T}_{n, Y}^{(r, s)}  \right) = \E \left( Y^{(u)}Y^{(r)} \right) \cdot A.
\end{align}
\end{lemma}

The following two corollaries follow immediately; for the sake of formality, their proofs are provided in Subsections \ref{proof_corollary_XY} and \ref{proof_corollary_XX} of Appendix \ref{Appendix_proofs}.

\begin{cor}\label{corollary_XY}
    Under MCAR, it holds that
    \begin{align}
    \lim_{n \to \infty} \Cov& \left( \sqrt{n} T_{n,X}^{(u,v)}, \sqrt{n} \hat{T}_{n,Y}^{(r,s)}  \right) \nonumber \\
    &= \Cov \left( X^{(u)}, Y^{(r)}  \right)  \bigg( \E \left( R^{(v)}R^{(r)} R^{(s)} \right) + \E \left( R^{(v)} \right)\E \left( R^{(r)} \right) \E \left( R^{(s)} \right) \nonumber \\
    & \hspace{150pt} - \E \left( R^{(v)} \right) \E \left( R^{(r)} R^{(s)} \right)   - \E \left( R^{(v)} R^{(r)} \right) \E \left( R^{(s)} \right) \bigg).  \label{cov_XY}
\end{align}
\end{cor}

\begin{cor}\label{corollary_XX}
    Under MCAR, it holds that
    \begin{align}\label{cov_XX}
    \lim_{n \to \infty} \Cov& \left( \sqrt{n} T_{n,X}^{(u,v)}, \sqrt{n} T_{n,X}^{(r,s)}  \right) = \Cov \left( X^{(u)}, X^{(r)}\right) \Cov \left( R^{(v)}, R^{(s)}\right).
    \end{align}
\end{cor}

The following result summarizes our findings. The proof, although quite obvious, is provided in Subsection \ref{proof_main} of Appendix \ref{Appendix_proofs}.

\begin{theorem}\label{theorem_main}
    Assume the data are MCAR, all variables have finite fourth moments, and \textit{either the incomplete variables have zero means or the response indicators are uncorrelated}. Then 
    \begin{align*}
    \sqrt{n}\bigg(   T_{n,X}^{(1,1)}, \dots,  &T_{n,X}^{(1, q)},  T_{n,X}^{(2,1)}, \dots,  T_{n,X}^{(2,q)}, \dots, T_{n,X}^{(p,1)}, \dots, T_{n,X}^{(p,q)},  \nonumber \\
    &  \hat{T}_{n, Y}^{(1,2)}, \dots, \hat{T}_{n, Y}^{(1,q)}, \hat{T}_{n, Y}^{(2,1)}, \hat{T}_{n, Y}^{(2,3)}, \dots, \hat{T}_{n, Y}^{(2,q)}, \dots, \hat{T}_{n, Y}^{(q,1)}, \dots, \hat{T}_{n, Y}^{(q, q-1)}  \bigg) \overset{D}{\to} \mathcal{N} (0, \Lambda)
\end{align*}
and
\begin{equation}\label{ultimate_convergence}
    A_n' := n\bigg(  T_{n,X}^{(1,1)}, \dots, \hat{T}_{n, Y}^{(q, q-1)}  \bigg) \hat{\Lambda}^{-1} \bigg(  T_{n,X}^{(1,1)}, \dots, \hat{T}_{n, Y}^{(q, q-1)}  \bigg)^T \overset{D}{\to} \chi^2_{pq + q(q-1)},
\end{equation}
where $\Lambda$ is corresponding limiting covariance matrix with limiting covariances \eqref{cov_YY_omitted}, \eqref{cov_XY}, and \eqref{cov_XX}, and $\hat{\Lambda}$ is its standard bias-adjusted estimate.
\end{theorem}

The convergence \eqref{ultimate_convergence} can be subsequently used to construct the rejection region and compute the $p$-value of the improved MCAR test.

\begin{remark}\label{remark_centering_problem}
    The assumption of either the incomplete variables having zero means or the response indicators being uncorrelated makes Theorem \ref{theorem_main} inapplicable in most practical scenarios. In real-world applications, we must mitigate this issue by centering the data before conducting the test. However, since we center the variables using available-case sample estimates of their means rather than the true, unknown theoretical ones, this introduces a potential methodological concern. Because these sample estimators are themselves random variables that converge at a $\sqrt{n}$-rate, plugging them into the test statistics injects an extra source of sampling variance into the system. 

    While simulations can evaluate how robust the unadjusted statistic $A_n'$ is to this plug-in effect, analyzing it theoretically provides a clearer understanding of the underlying mechanics. Keeping the baseline derivations of Lemma \ref{lemma_YY} and Theorem \ref{theorem_main} is useful because they show how the data and the missingness mechanism interact under ideal conditions. 
\end{remark}

To handle practical scenarios where the data must be centered using sample means, we use a first-order Taylor expansion (or Hoeffding projection) to account for the exact variance introduced by this estimation. This allows us to drop the restrictive assumptions of Theorem \ref{theorem_main} and establish a fully generalized framework. The proof can be found in Subsection \ref{proof_theorem_plugin} of Appendix \ref{Appendix_proofs}.

\begin{theorem}\label{theorem_plugin}
Let $\boldsymbol{\hat{\mu}}_Y = (\hat{\mu}_Y^{(1)}, \dots, \hat{\mu}_Y^{(q)})^T$ be the available-case sample means defined by $\hat{\mu}_Y^{(u)} = \sum_{i=1}^n Y_i^{(u)} R_i^{(u)}/\hat{n}^{(u)}$, for $u = 1,2, \dots, q$. Let $\hat{T}_{n, Y}^{(u, v)}(\boldsymbol{\hat{\mu}}_Y)$ denote the test statistics computed on the sample-centered variables:
\begin{align*}
    \hat{T}_{n, Y}^{(u, v)}(\boldsymbol{\hat{\mu}}_Y) &= \frac{1}{n} \sum_{i=1}^n \left(Y_i^{(u)} - \hat{\mu}_Y^{(u)}\right) R_i^{(u)} R_i^{(v)} - \frac{1}{n(n-1)} \sum_{i=1}^n \sum_{\substack{j=1 \\ j \neq i}}^n \left(Y_i^{(u)} - \hat{\mu}_Y^{(u)}\right) R_i^{(u)} R_j^{(v)}.
\end{align*}
Under the null hypothesis of MCAR and the assumption of finite fourth moments, as $n \to \infty$, we have
\begin{align*}
    \sqrt{n}\bigg(   T_{n,X}^{(1,1)}, \dots,  &T_{n,X}^{(p,q)},  \hat{T}_{n, Y}^{(1,2)}(\boldsymbol{\hat{\mu}}_Y), \dots, \hat{T}_{n, Y}^{(q, q-1)}(\boldsymbol{\hat{\mu}}_Y)  \bigg) \overset{D}{\to} \mathcal{N} (0, \Omega),
\end{align*}
where the elements of the limiting covariance matrix $\Omega$ are given by
\begin{align}\label{omega_generic}
    \Omega_{(u,v), (r,s)} = \Cov\left(Z^{(u)}, Z^{(r)}\right) \cdot \E\left[ R^{(u)} R^{(r)} \left( R^{(v)} - \frac{\E(R^{(u)}R^{(v)})}{\E(R^{(u)})} \right) \left( R^{(s)} - \frac{\E(R^{(r)}R^{(s)})}{\E(R^{(r)})} \right) \right],
\end{align}
where $Z^{(u)}$ represents either a fully observed variable $X^{(u)}$ (with $R^{(u)} \equiv 1$) or a partially observed variable $Y^{(u)}$. Consequently, the plug-in test statistic
\begin{equation}\label{ultimate_convergence_plugin}
    A_n'' := n\bigg(  T_{n,X}^{(1,1)}, \dots, \hat{T}_{n, Y}^{(q, q-1)}(\boldsymbol{\hat{\mu}}_Y)  \bigg) \hat{\Omega}^{-1} \bigg(  T_{n,X}^{(1,1)}, \dots, \hat{T}_{n, Y}^{(q, q-1)}(\boldsymbol{\hat{\mu}}_Y)  \bigg)^T \overset{D}{\to} \chi^2_{pq + q(q-1)}.
\end{equation}
\end{theorem}

\begin{remark}
Equation \eqref{omega_generic} provides an elegant and unified representation of the asymptotic covariance. Notice that if $Z^{(u)} = X^{(u)}$, then $R^{(u)} \equiv 1$, and the inner fractional term reduces cleanly to $\E(R^{(v)})$, which perfectly matches the covariance structure derived under the simpler conditions of Corollary \ref{corollary_XY} and Corollary \ref{corollary_XX}.
\end{remark}

\begin{remark}\label{remark_invariance}
    For any fully observed variable $X^{(u)}$ (where $R^{(u)} \equiv 1$), the corresponding test statistic $T_{n,X}^{(u,v)}$ is strictly invariant to any constant location shift. If we replace $X_i^{(u)}$ with $X_i^{(u)} - C$ for an arbitrary constant $C$, the shift cancels out completely:
    \begin{align*}
        \frac{1}{n} \sum_{i=1}^n (X_i^{(u)} - C)R_i^{(v)} &- \frac{1}{n(n-1)} \sum_{i=1}^n \sum_{\substack{j=1 \\ j \neq i}}^n (X_i^{(u)} - C)R_j^{(v)} \\
        &= T_{n,X}^{(u,v)} - \frac{C}{n}\sum_{i=1}^n R_i^{(v)} + \frac{C}{n(n-1)} (n-1) \sum_{j=1}^n R_j^{(v)} \\
        &= T_{n,X}^{(u,v)} - \frac{C}{n}\sum_{i=1}^n R_i^{(v)} + \frac{C}{n}\sum_{j=1}^n R_j^{(v)} = T_{n,X}^{(u,v)}.
    \end{align*}
    Consequently, while centering the incomplete variables $Y^{(u)}$ by $\hat{\mu}_Y^{(u)}$ alters their asymptotic influence functions and requires the modified covariance matrix $\Omega$, centering the complete variables $X^{(u)}$ by their sample means has absolutely no effect on their values. Doing so could sometimes be an operational choice that allows us to treat all variables uniformly in the computational implementation.
\end{remark}

\section{Empirical study}\label{sec:simulation}

In this section, we present the results of an extensive simulation study which is conducted to examine how the novel test behaves in terms of empirical type I error and power. The novel test is compared to Little's MCAR test, as well as the old test based on the statistic $A_n$ from \eqref{An}, which novel $A_n''$ from \eqref{ultimate_convergence} improves upon. 

As previously noted, the original test based on $A_n$ could only detect correlations between the response indicators and fully observed variables. Therefore, in addition to evaluating alternatives undetectable by the original test, we also compare $A_n$, Little's $d^2$, and the proposed $A_n''$ in these scenarios \citep{aleksic2024novel}. This allows us to assess whether the improvements offered by the new test come with potential trade-offs, such as reduced power or calibration issues.

The missingness probability for any value in the data, i.e. the probability that a specific data cell is missing, ranges from $3\%$ to $30\%$. Furthermore, we have decided to study sample sizes of $n = 100$, $n = 200$, and $n = 300$, which appears to be adequate for illustrating the quality of asymptotic approximations.

Throughout the rest of this section, we use the abbreviation $iXjY$ to denote the dataset with a total of $i+j$ variables, where $i$ of them are complete, and $j$ of them are incomplete. We mostly present only the $2X3Y$ results in the main text to avoid overload; findings for other cases are consistent and thus omitted. They can be found in Supplementary material.

All simulations are performed with $N = 2000$ replications and at nominal level of $\alpha = 0.05$.

\subsection{Study design}

\subsubsection{Generating the data}

For the data distributions, we use the standard normal distribution, as well as the normal distribution with marginal means equal to 1 and covariance matrix $0.5I + 0.5J$, where $I$ is the identity matrix and $J$ is a matrix with all elements equal to $1$. We also consider a Clayton copula with parameter 1 and exponential $\mathcal{E}(1)$ margins, as well as $\chi_4^2$ margins \citep[see, e.g.,][]{fischer2012constructing}. The main idea behind this choice is that Little's test relies on the normality assumption. Given this, it is important to assess the performance of the novel test for normally distributed data, for data whose distribution deviates substantially from normality, and intermediate cases. The Clayton copula was also used in an independence testing scenario by \cite{cuparic2024ipcw}, where the test of independence by \cite{kochar1990distribution} was adapted for the setting of randomly censored data.

\subsubsection{Generating missingness with uncorrelated response indicators}

For implementation, \texttt{R} package \texttt{missMethods} is used \citep{missMethods}. For the null distribution case, function \texttt{delete\_MCAR} is used. We stick to the alternatives implemented in functions \texttt{delete\_MAR\_1\_to\_x}, with recommended choice of $x = 9$ and argument {\texttt{n\_mis\_stochastic = FALSE}}, and \texttt{delete\_MAR\_rank}. The main idea between these mechanisms is that, for each incomplete data column, we have the so-called \emph{control column}, that is fully observed, and the data from that column is used to dictate the missingness probability in the incomplete one. {The first} mechanism works by setting a specific threshold (default is median), and then splitting the cases into two groups: those that have value of control variable smaller than the threshold, and those that do not. Then, the missingness is introduced such that the odds of a value being missing in those two groups are $1:x$. For the second mechanism, the probability that a value is missing is directly proportional on the rank of its observed pair. For more details, {we} refer to \cite{santos2019generating}, where they were first introduced.

\subsubsection{Generating missingness with correlated response indicators}

{To examine the behavior of the improved test when centering is required, we generate missing data with correlated response indicators, controlling their correlation coefficient. This is done by modifying the functions \texttt{delete\_MCAR}, \texttt{delete\_MAR\_1\_to\_x}, and \texttt{delete\_MAR\_rank}. Algorithm \ref{algorithm_missingness} illustrates the procedure for the $2X3Y$ data and positive correlation, and it can be easily adapted to other dimensions and correlation structures.

\begin{algorithm}
  \caption{Generating missingness with positively correlated response indicators.}\label{algorithm_missingness}
  \begin{algorithmic}[1]
    \State Start with the complete sample $(x_1, x_2,  \dots, x_n)$, where each $x_j = \left(x_j^{(1)}, x_j^{(2)}, y_j^{(1)}, y_j^{(2)}, y_j^{(3)}\right)$;
    \State Specify the desired missingness probability $p$ and {the} correlation coefficient {$\rho$};
    \State Generate missingness in variables $y^{(1)}$ and $y^{(2)}$ using {the} probability $p$ and a chosen method; record the response indicator vector $r^{(2)}$;
    \State Generate {a} random vector $r$ of length $n$ {consisting} of zeros and ones, where $0$ appears with probability $p$, and $1$ with probability $1 - p$;
    \State Generate {a} response indicator vector $r^{(3)}$ of length $n$ such that its $j$th element is equal to $r^{(2)}_j$ with probability $r$, and $r_j$ with probability $1-r$;
    \State Generate missingness in $y^{(3)}$ according to the response indicator $r^{(3)}$.
  \end{algorithmic}
\end{algorithm}

The following lemma formally establishes that the correlated response indicators generated in Algorithm \ref{algorithm_missingness} have a correlation coefficient equal to {$\rho$}. The proof can be found in Subsection \ref{appendix_lemma_algorithm} of Appendix \ref{Appendix_proofs}.
\begin{lemma}\label{lemma_algorithm}
    Let $R_1$ and $R_2$ be two independent indicator random variables with the same expected value $q$. Let $U$ be a random variable uniformly distributed on $[0,1]$ and independent of both $R_1$ and $R_2$. If, for {$\rho \in [0,1]$}, 
    \begin{align*}
        R_3 = I\left\{ U \leq  {\rho}\right\} R_1 + I\left\{ U > {\rho} \right\}R_2, 
    \end{align*}
    then $\mathbf{Cor} (R_1, R_3) = {\rho}$.
\end{lemma}

\begin{remark}\label{remark_alg_mis}
    Similarly to the missingness settings shown in Figures \ref{fig:2X3Y_stdnorm_MAR_rank_undetectable1}–\ref{fig:2X3Y_clayton1_chisq4_MAR_rank_undetectable1}, we adapt Algorithm \ref{algorithm_missingness} so that the variable $Y^{(1)}$ governs the missingness of $Y^{(2)}$, while the missingness in $Y^{(3)}$ is generated to be correlated with that of $Y^{(2)}$. Additionally, MCAR is imposed on $Y^{(1)}$ to make the alternative hardly detectable for the $A_n$-based test.
\end{remark}

\subsection{Performance of tests under zero mean or uncorrelated response indicators}

In this subsection, we present the results of simulations in which the $A_n''$-based testing is conducted under the assumption that the variables have zero means or that the response indicators are uncorrelated. In this case, centering the data  prior to testing is unnecessary for its validity; however, to have a unified framework, data are centered nevertheless. As seen from Remark \ref{remark_invariance}, it has no impact on the procedure. We first present results for the exact missingness scenarios considered in \cite{aleksic2024novel}, and subsequently for scenarios in which the alternative was undetectable or only marginally detectable by the original test.

\subsubsection{Performance in scenarios where $A_n$-based test was compared to Little's test}\label{subsec:simulation_old}

\begin{figure}
    \centering
    \includegraphics[width=\textwidth]{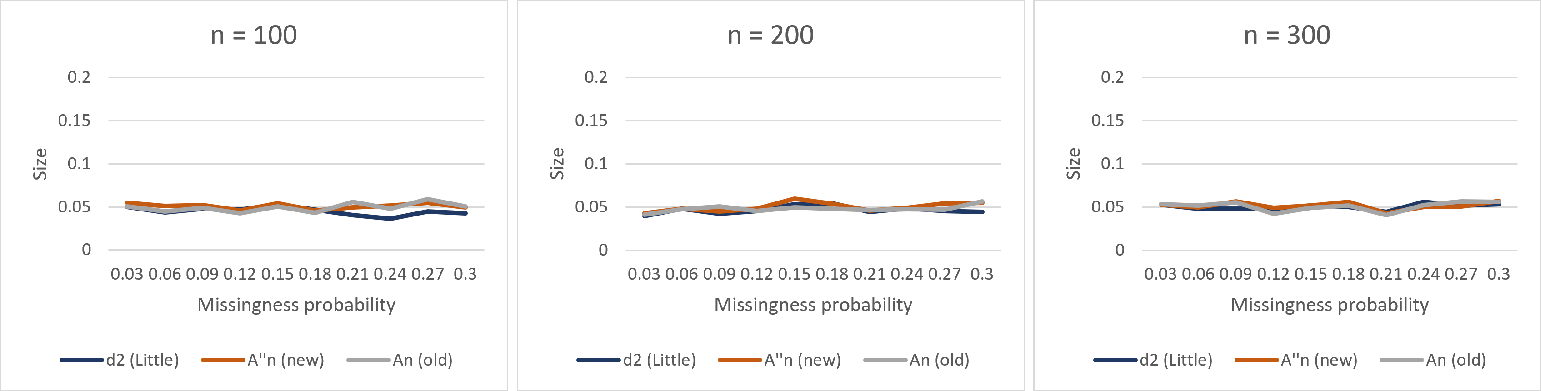}
    \caption{Empirical type I errors for 2X3Y case, standard normal distribution.}
    \label{fig:2X3Y_stdnorm_MCAR}
\end{figure}

\begin{figure}
    \centering
    \includegraphics[width=\textwidth]{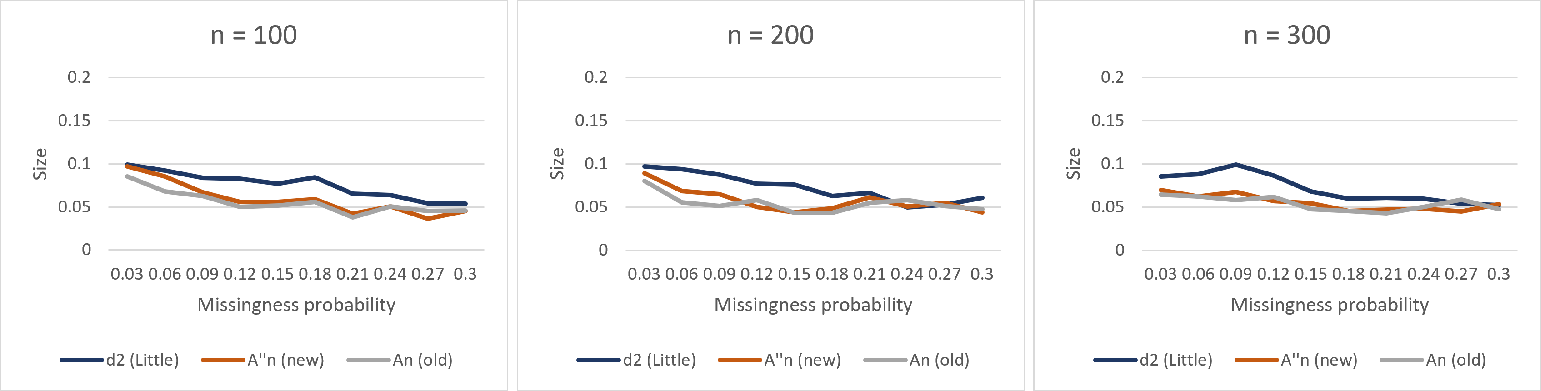}
    \caption{Empirical type I errors for 2X3Y case, Clayton copula with parameter 1 and $\mathcal{E}(1)$ margins.}
    \label{fig:2X3Y_clayton1_exp1_MCAR}
\end{figure}

\begin{figure}
    \centering
    \includegraphics[width=\textwidth]{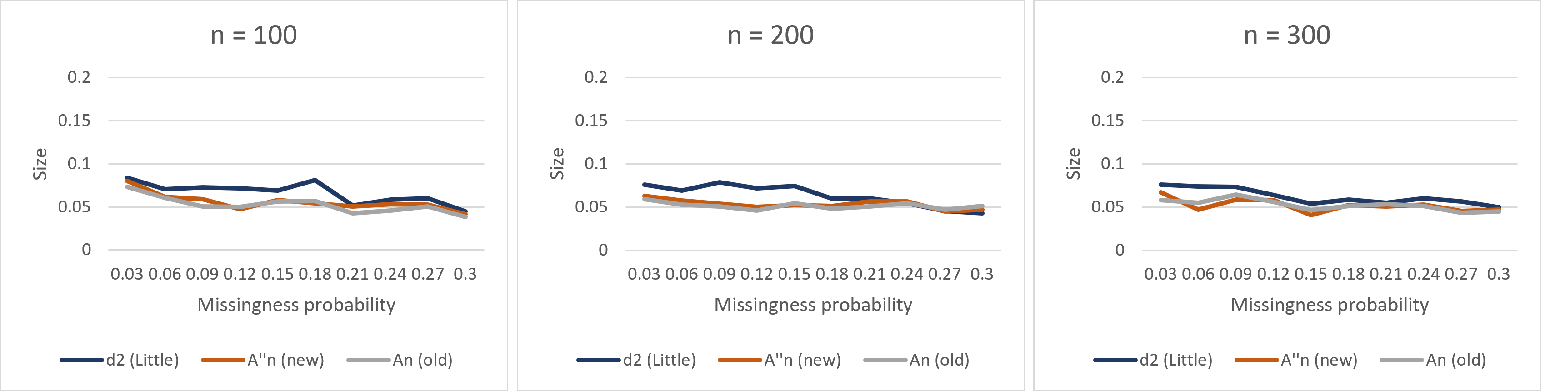}
    \caption{Empirical type I errors for 2X3Y case, Clayton copula with parameter 1 and $\chi^2_4$ margins.}
    \label{fig:2X3Y_clayton1_chisq4_MCAR}
\end{figure}

\begin{figure}
    \centering
    \includegraphics[width=\textwidth]{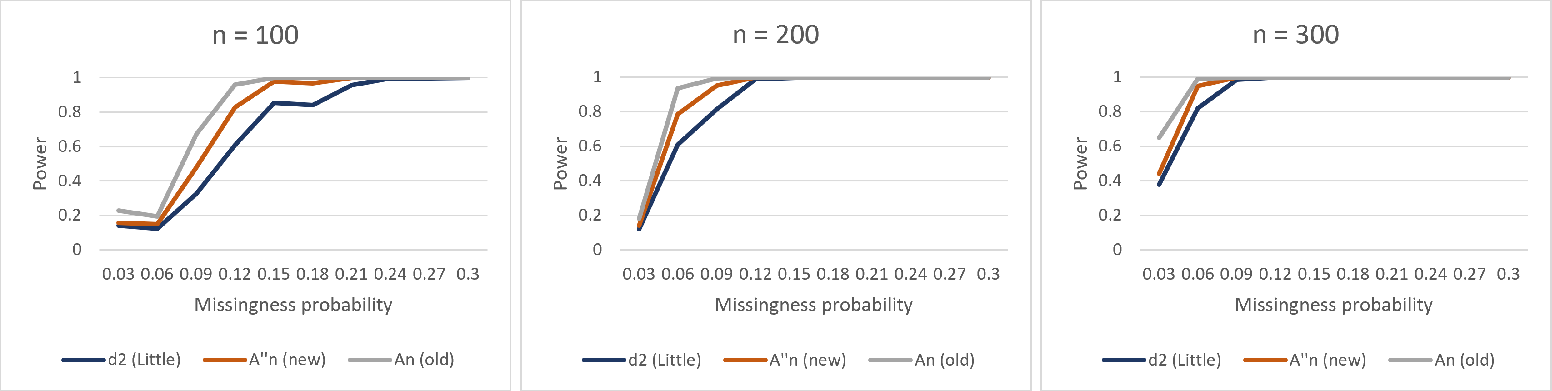}
    \caption{Empirical test powers for 2X3Y case, standard normal distribution, MAR 1 to 9 (var. 1 controls missingness in var. 3 and 
 var. 5, var. 2 controls var. 4).}
    \label{fig:2X3Y_stdnorm_MAR_1_to_9}
\end{figure}

\begin{figure}
    \centering
    \includegraphics[width=\textwidth]{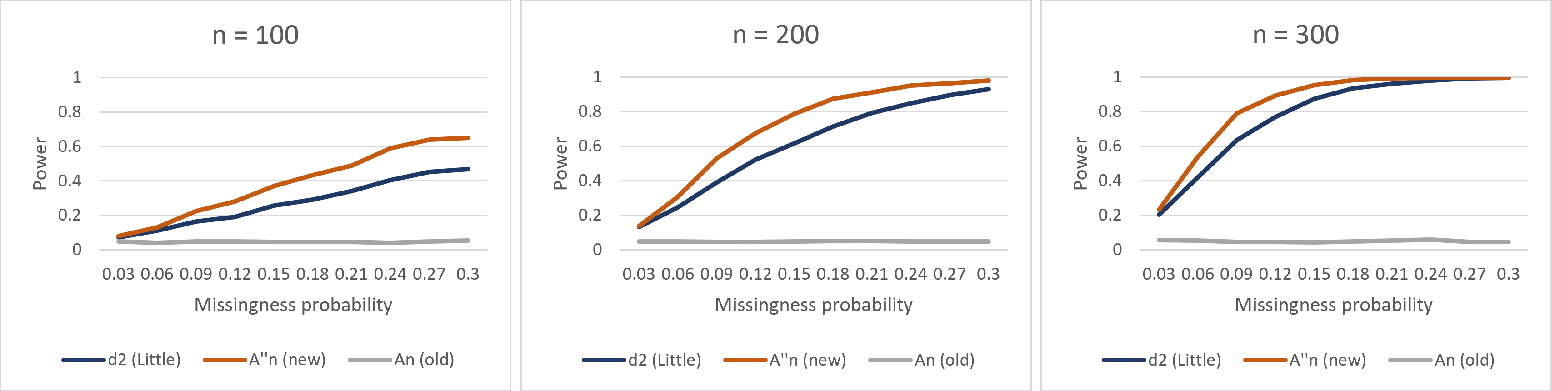}
    \caption{Empirical test powers for 2X3Y case, standard normal distribution, combination of MAR rank and MCAR (var. 3 controls missingness in var. 4 and 
 var. 5, and then MCAR missingness is generated in var. 3).}
    \label{fig:2X3Y_stdnorm_MAR_rank_undetectable1}
\end{figure}

As demonstrated in Figure \ref{fig:2X3Y_stdnorm_MCAR}, all three tests are well-calibrated under the standard normal distribution, with empirical type I error rates closely matching the nominal level. However, Figures \ref{fig:2X3Y_clayton1_exp1_MCAR} and \ref{fig:2X3Y_clayton1_chisq4_MCAR} reveal that Little's $d^2$ exhibits type I error inflation. This deviation is particularly pronounced for exponentially distributed margins, where the error rate reaches nearly twice the nominal level, a behavior most notable at small-to-moderate missingness rates (Figure \ref{fig:2X3Y_clayton1_exp1_MCAR}). Conversely, $A_n$ and $A_n''$ display very similar performance and remain much better calibrated than $d^2$. This advantage is most apparent in Figure \ref{fig:2X3Y_clayton1_exp1_MCAR}, where the underlying data distribution deviates most severely from normality.

Under the \textit{MAR 1 to x} alternative with normal data, Figure \ref{fig:2X3Y_stdnorm_MAR_1_to_9} indicates that the novel test based on $A_n''$ experiences a loss of power compared to the original one based on $A_n$. Nevertheless, it consistently outperforms Little's MCAR test, especially for smaller sample sizes. Similar trends are observed across the other underlying distributions, as well as under the MAR rank mechanism. To maintain conciseness, these results are presented in the Supplementary material.

\subsubsection{Performance in novel scenarios}\label{subsec:simulation_novel}

\begin{figure}
	\centering
	\includegraphics[width=\textwidth]{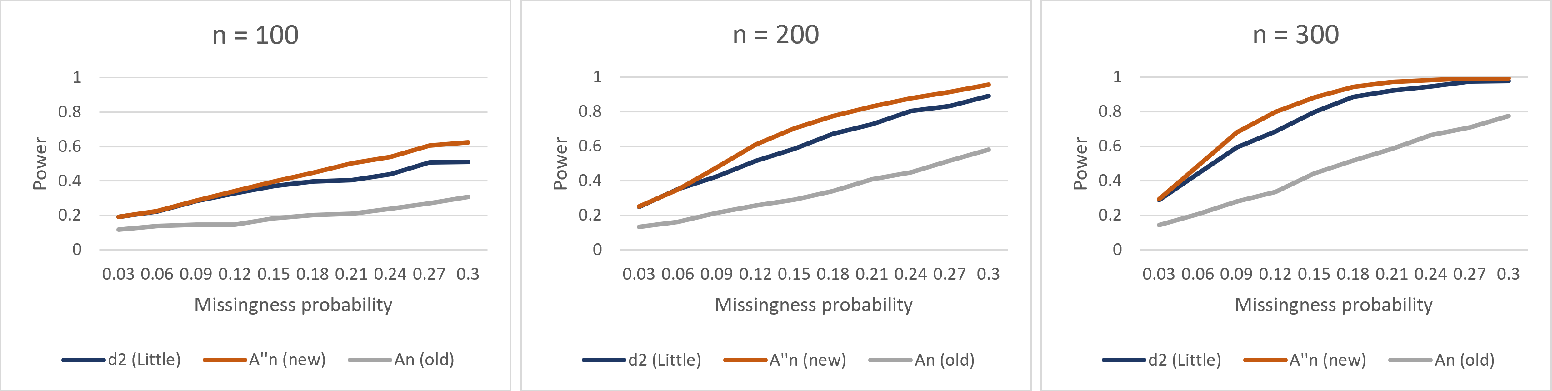}
	\caption{Empirical test powers for 2X3Y case, Clayton copula with parameter 1 and $\mathcal{E}(1)$ margins, combination of MAR rank and MCAR (var. 3 controls missingness in var. 4 and 
		var. 5, and then MCAR missingness is generated in var. 3).}
	\label{fig:2X3Y_clayton1_exp1_MAR_rank_undetectable1}
\end{figure}

\begin{figure}
	\centering
	\includegraphics[width=\textwidth]{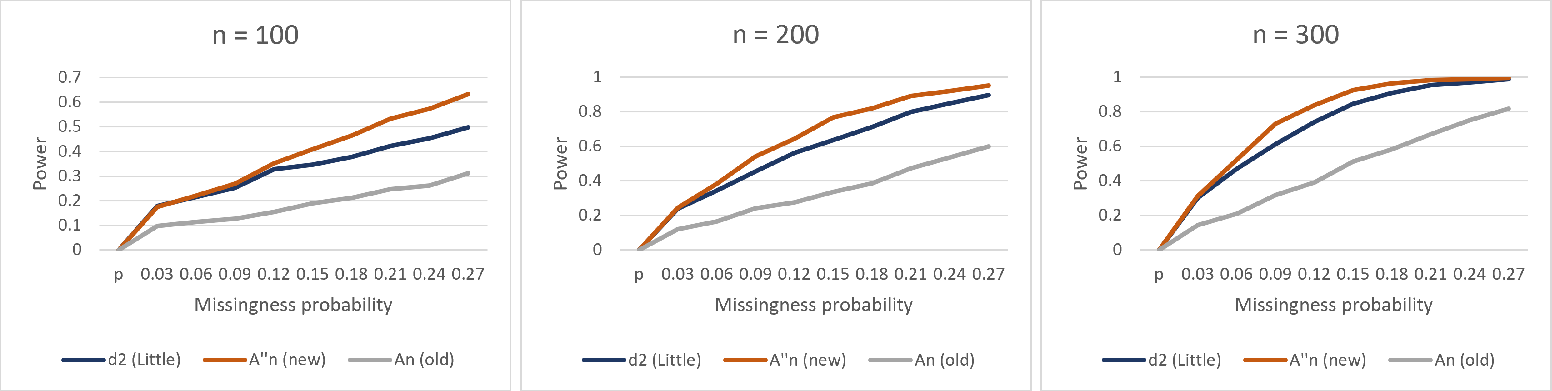}
	\caption{Empirical test powers for 2X3Y case, Clayton copula with parameter 1 and $\chi^2_4$ margins, combination of MAR rank and MCAR (var. 3 controls missingness in var. 4 and 
		var. 5, and then MCAR missingness is generated in var. 3).}
	\label{fig:2X3Y_clayton1_chisq4_MAR_rank_undetectable1}
\end{figure}

Figure \ref{fig:2X3Y_stdnorm_MAR_rank_undetectable1} shows the power performance for a specific MAR setting using standard normal $2X3Y$ data: variable 3 controls the missingness in variables 4 and 5 according to the \textit{MAR rank} mechanism, and MCAR missingness is subsequently generated in variable 3. This is a representative example of a setting where response indicators depend on an incomplete column, which constitutes an alternative that is undetectable by the old test. The novel test once again outperforms Little's test. 

Figure \ref{fig:2X3Y_clayton1_exp1_MAR_rank_undetectable1} demonstrates that the old test can detect the alternative when the variables are correlated, specifically for the Clayton copula with parameter 1 and $\mathcal{E}(1)$ margins; the old test captures the dependence through the completely observed variables. However, the old test exhibits significantly lower power, whereas the novel test is slightly more powerful than Little's test. For $\chi^2_4$ margins (Figure \ref{fig:2X3Y_clayton1_chisq4_MAR_rank_undetectable1}), the results are similar. The behavior observed in Figures \ref{fig:2X3Y_stdnorm_MAR_rank_undetectable1},  \ref{fig:2X3Y_clayton1_exp1_MAR_rank_undetectable1}, and \ref{fig:2X3Y_clayton1_chisq4_MAR_rank_undetectable1} persists across different dimensions, distributions, and alternatives that are undetectable by the $A_n$-based test. Further examples of this behavior are provided in the Supplementary material.

\begin{remark}
	It is important to note that in scenarios where the data deviate from normality, Little’s test tends to reject the null hypothesis more frequently than it should, indicating inflated type I error rates. Consequently, the apparent power of Little’s test under non-normal settings should be interpreted with caution, as part of the observed rejections may come from this liberal behavior rather than genuine departures from the null. In other words, the reported powers for non-normal data likely overestimate the true power of the test. Therefore, when comparing the performance of the proposed methods, the results obtained under normality should be regarded as the most reliable benchmark, providing the best estimate of the actual differences in power between procedures.
\end{remark}

Since all three studied tests {rely on the assumption of all of the variables having finite fourth moments}, it is interesting to examine the robustness of tests when that assumption is not fulfilled. Figure \ref{fig:2X3Y_stdt2_MCAR} shows the empirical type I errors for the standard Student's $t$-distribution with 2 degrees of freedom, which does not have finite fourth moments. As we can see, the novel test performs much better that Little's test, even for larger sample sizes. Despite the tendency of $d^2$-based test to reject the null hypothesis in this setting {even when the null hypothesis is true}, Figure \ref{fig:2X3Y_stdt2_MAR_1_to_9} shows that the novel test is significantly more powerful. 

\begin{figure}
	\centering
	\includegraphics[width=\textwidth]{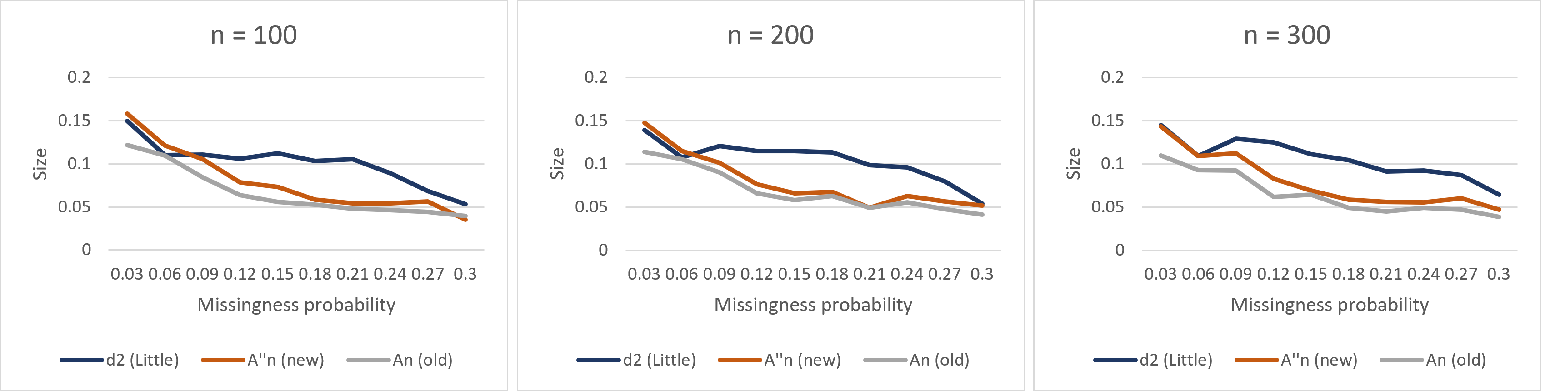}
	\caption{Empirical test sizes for 2X3Y case, standard Student's $t_2$ distribution}
	\label{fig:2X3Y_stdt2_MCAR}
\end{figure}

\begin{figure}
	\centering
	\includegraphics[width=\textwidth]{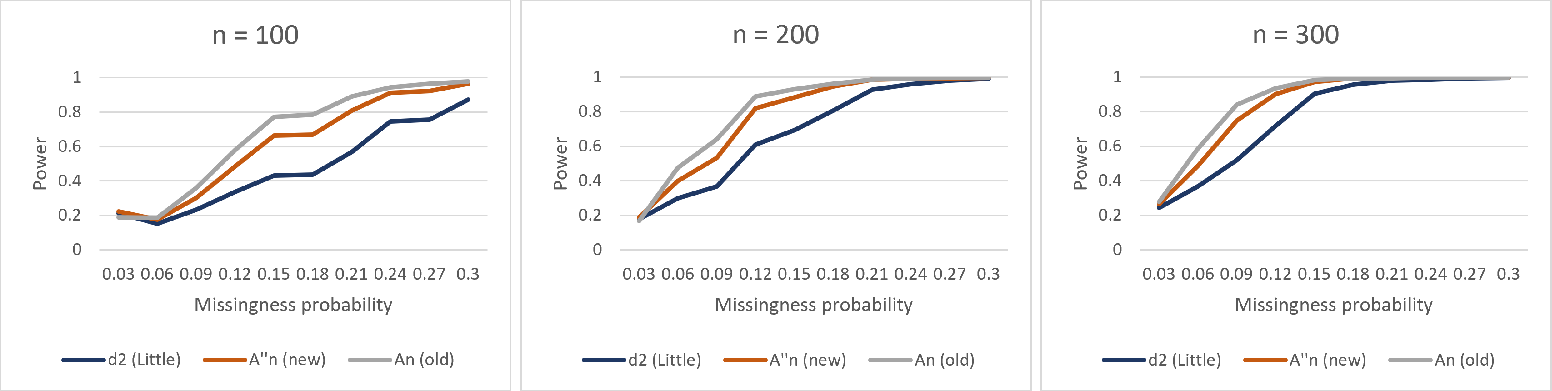}
	\caption{Empirical test powers for 2X3Y case, standard Student's $t_2$ distribution, MAR 1 to 9 (var. 1 controls missingness in var. 3 and 
		var. 5, var. 2 controls var. 4)}
	\label{fig:2X3Y_stdt2_MAR_1_to_9}
\end{figure}

\begin{figure}
	\centering
	\includegraphics[width=\textwidth]{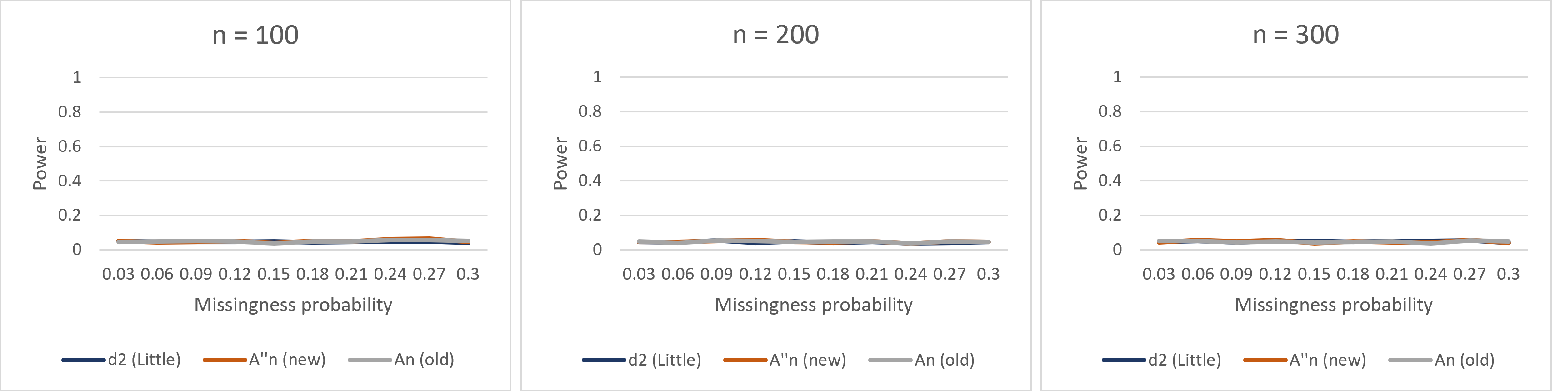}
	\caption{Empirical test powers for 2X3Y case, standard normal distribution, MNAR (upper) censoring}
	\label{fig:2X3Y_stdnorm_MNAR_upper_censoring}
\end{figure}

\begin{figure}
	\centering
	\includegraphics[width=\textwidth]{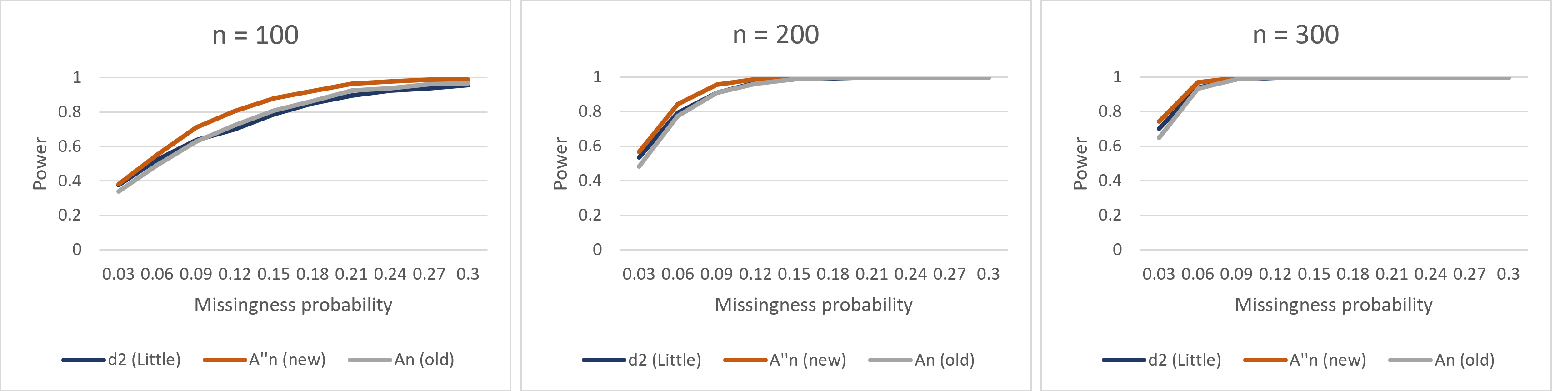}
	\caption{Empirical test powers for 2X3Y case, Clayton copula with parameter 1 and $\mathcal{E}(1)$ margins, MNAR (upper) censoring}
	\label{fig:2X3Y_clayton1_exp1_MNAR_upper_censoring}
\end{figure}

\begin{figure}
	\centering
	\includegraphics[width=\textwidth]{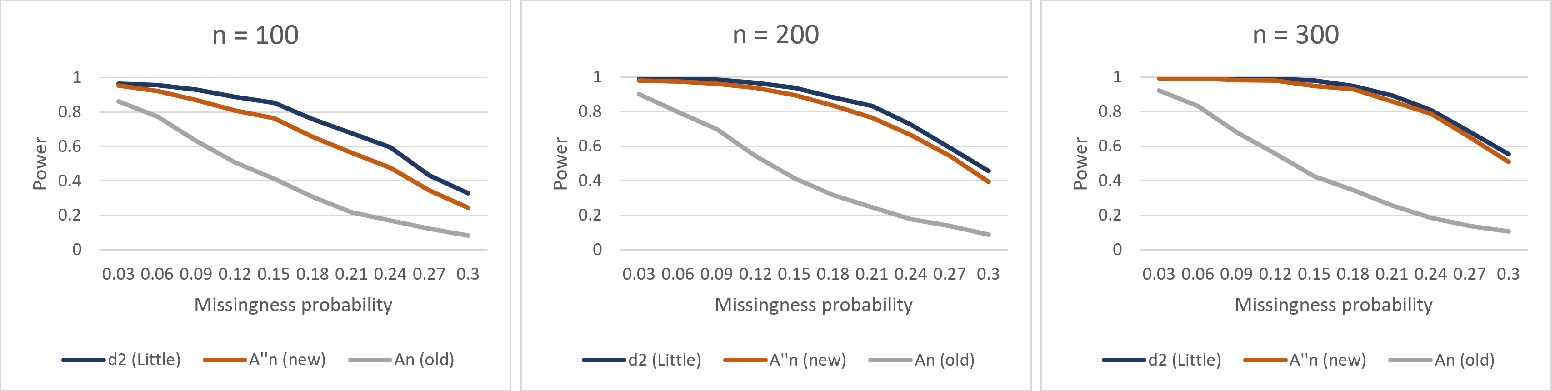}
	\caption{Empirical test powers for 2X3Y case, standard Student's $t_2$ distribution, MNAR (upper) censoring}
	\label{fig:2X3Y_stdt2_MNAR_upper_censoring}
\end{figure}

Another important scenario in which the novel test needs to be examined is the case of MNAR data. Since the test by its construction is not able to calculate {the} covariance between the incomplete variable and its response indicators, alternatives that are "purely MNAR" should be undetectable for the test. More precisely, those are alternatives where the only form of dependence between the response indicators and the data is realized between the variable and its indicators, but not any others. However, as seen from Figure \ref{fig:2X3Y_stdnorm_MNAR_upper_censoring}, in {the} case of {the} standard normal distribution, Little's test is not able to detect such alternative either. Figure \ref{fig:2X3Y_clayton1_exp1_MNAR_upper_censoring} presents behavior in the case of same missingness mechanism, but Clayton copula with parameter 1 and exponential {with parameter 1} margins. As we can see, all tests have practically the same power, and are able to detect the alternative. The same behavior is noted for other studied distributions and dimensions.

However, there are exceptions that behave unexpectedly, such as previously studied Student's $t_2$ distribution which does not have finite fourth moments. When combined with upper censoring as the only missingness mechanism, Figure \ref{fig:2X3Y_stdt2_MNAR_upper_censoring} {shows} that all three test experience \textit{loss} of power {as the missingness rate increases}, which is not expected, and was not observed in the scenarios before. {It appears that the combination of an undetectable alternative and data from a population that do not satisfy the assumption of finite fourth moment is too challenging for the tests to handle}, and they start behaving in an unexpected manner. We have tried replacing the identity scale matrix of the standard $t_2$ distribution with the matrix that has unit diagonal elements, and others equal to 0.1 and 0.5, respectively, 
but it did not help the Little's test, and behavior persisted. For example, for the scale matrix with non-diagonal elements equal to 0.5, the novel test stopped having decreasing power for $n = 100$, but Little's test stabilized for $n = 300$.  These results are part of the Supplementary material.

\begin{remark}\label{remark_transform}
	The novel test we have presented is based on estimating the covariance between response indicators and data variables, which is a measure of linear dependence. To capture other forms of dependence, one is free to transform the variables and apply the test to the transformed data, if some other form of dependence is expected to occur based on previous experience. One such example where transforming the variables improves the power performance of the test can be found in \cite{aleksic2024novel}, where the original $A_n$-based test was introduced. The transformation of variables in this context could be the best solution in those cases where specific dependence between the data and the response indicators is to be expected. 
\end{remark}

\subsection{Performance of tests under nonzero mean and correlated response indicators}

This subsection examines the properties of the $A_n''$-based test in situations where data centering is required, that is, when some variables have nonzero means and certain response indicators are correlated. We provide a concise overview of the conducted simulations to avoid overloading the text. For this purpose, only some of the results for the normal distribution with all means equal to $1$ and covariance matrix $0.5I + 0.5J$ are presented. The behavior observed for other distributions is consistent with the results shown here, except that Little’s test exhibits weaker type I error control for distributions that deviate substantially from normality. This was also the case for the setting of zero means and uncorrelated response indicators.  We also note that the same pattern is also observed across different dimensions and alternative hypotheses. Those additional results can be found in the Supplementary material.

\begin{figure}
	\centering
	\includegraphics[width=\textwidth]{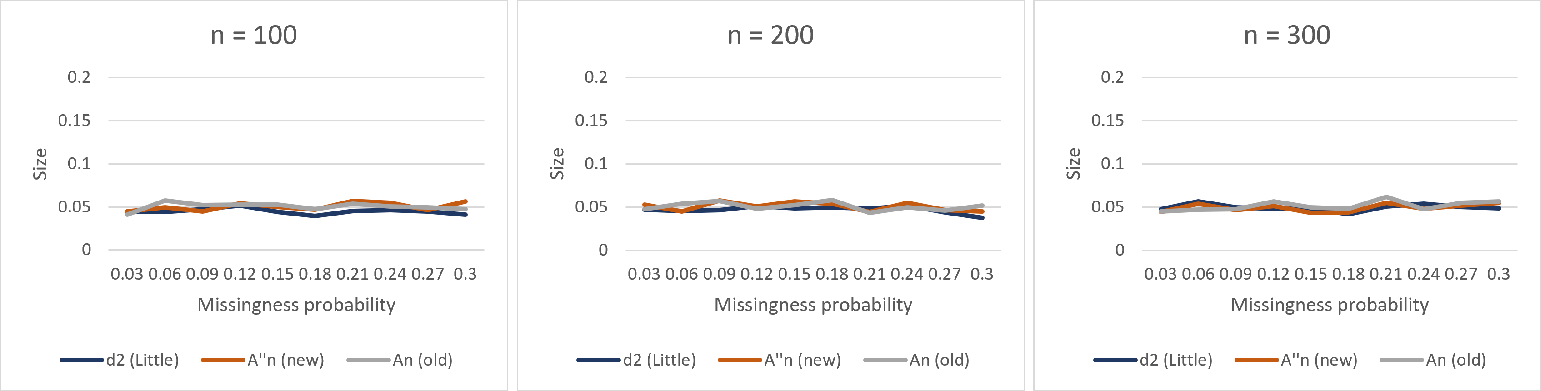}
	\caption{Empirical type I errors for 2X3Y case, normal distribution with mean $(1,1,1,1,1)$ and covariance matrix $0.5I_5+0.5J_5$, MCAR using Algorithm \ref{algorithm_missingness} with ${\rho}=0.2$}
	\label{fig:2X3Y_norm_mean1_Sigma05I05J_MCAR_r02}
\end{figure}

\begin{figure}
	\centering
	\includegraphics[width=\textwidth]{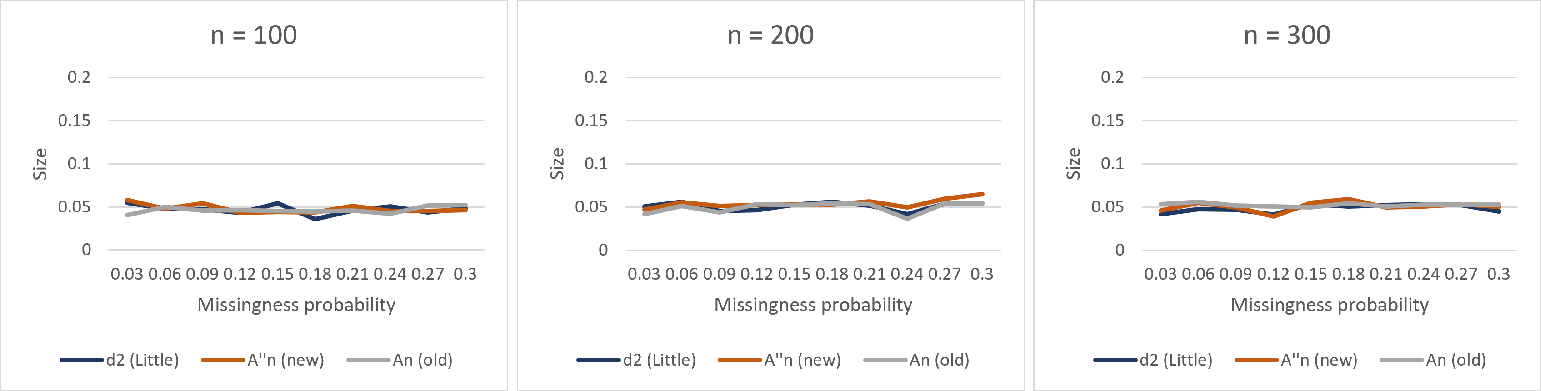}
	\caption{Empirical type I errors for 2X3Y case, normal distribution with mean $(1,1,1,1,1)$ and covariance matrix $0.5I_5+0.5J_5$, MCAR using Algorithm \ref{algorithm_missingness} with ${\rho}=0.8$}
	\label{fig:2X3Y_norm_mean1_Sigma05I05J_MCAR_r08}
\end{figure}

\begin{figure}
	\centering
	\includegraphics[width=\textwidth]{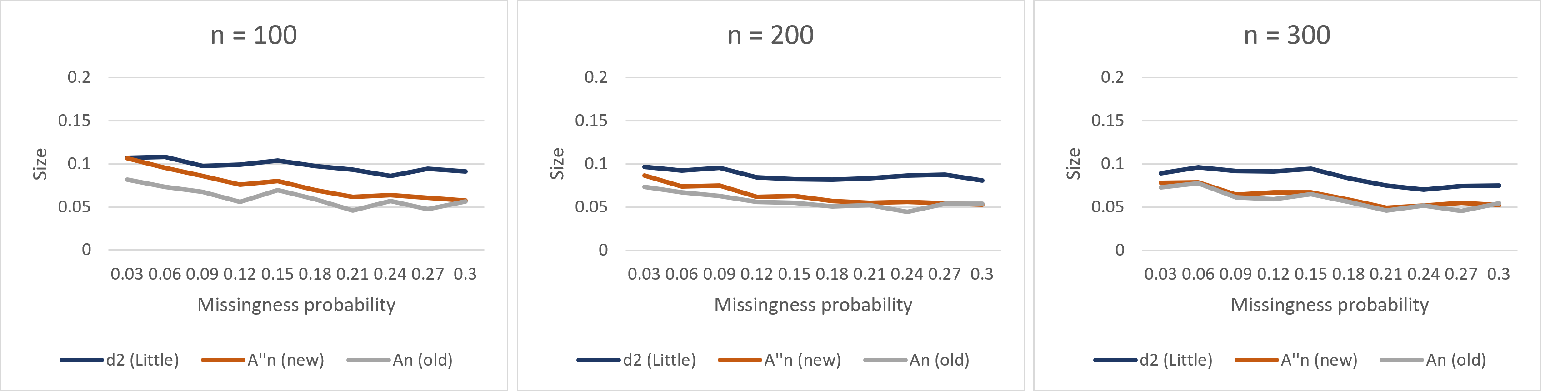}
	\caption{Empirical type I errors for 2X3Y case, Clayton copula with parameter $1$ and $\mathcal{E}(1)$ margins, MCAR using Algorithm \ref{algorithm_missingness} with ${\rho}=0.8$}
	\label{fig:2X3Y_clayton1_exp1_MCAR_r08}
\end{figure}

Figures \ref{fig:2X3Y_norm_mean1_Sigma05I05J_MCAR_r02} and \ref{fig:2X3Y_norm_mean1_Sigma05I05J_MCAR_r08} display the empirical type I errors for normally distributed data with correlated response indicators and nonzero means. As observed, all tests are well calibrated. As previously noted, Little’s MCAR test exhibits weaker control of the type I error when the data deviate from normality. Figure \ref{fig:2X3Y_clayton1_exp1_MCAR_r08} provides an illustrative example of such behavior.

\begin{figure}
	\centering
	\includegraphics[width=\textwidth]{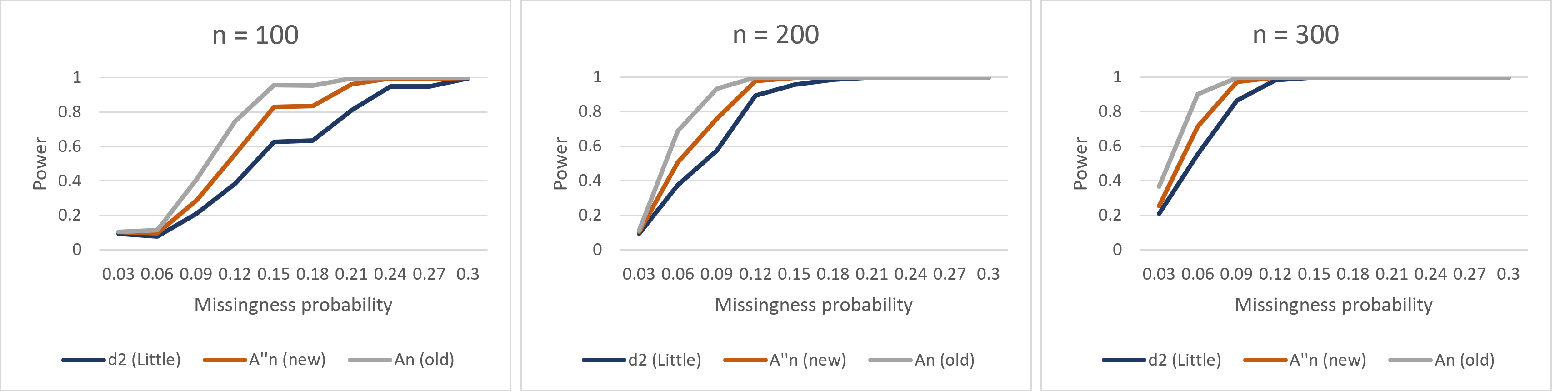}
	\caption{Empirical test powers for 2X3Y case, normal distribution with mean $(1,1,1,1,1)$ and covariance matrix $0.5I_5+0.5J_5$, MAR 1 to 9 using Algorithm \ref{algorithm_missingness} with ${\rho}=0.8$}
	\label{fig:2X3Y_norm_mean1_Sigma05I05J_MAR_1_to_9_r08}
\end{figure}

\begin{figure}
	\centering
	\includegraphics[width=\textwidth]{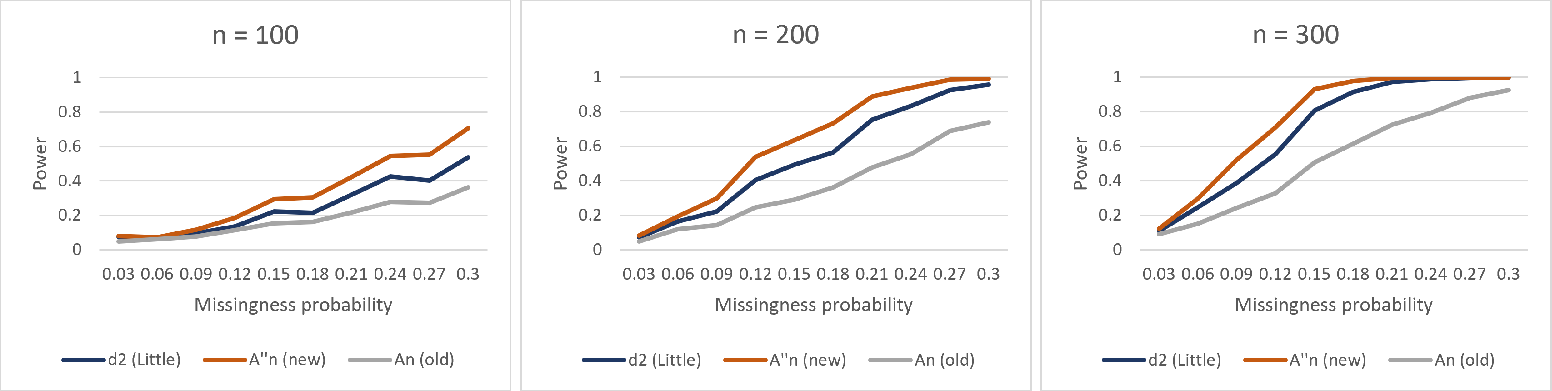}
	\caption{Empirical test powers for 2X3Y case, normal distribution with mean $(1,1,1,1,1)$ and covariance matrix $0.5I_5+0.5J_5$, MAR 1 to 9 using modified Algorithm \ref{algorithm_missingness} from Remark \ref{remark_alg_mis}, with ${\rho}=0.8$}
	\label{fig:2X3Y_norm_mean1_Sigma05I05J_MAR_1_to_9_masked_r08}
\end{figure}

Figure \ref{fig:2X3Y_norm_mean1_Sigma05I05J_MAR_1_to_9_r08} provides a representative example of the power behavior under an alternative detectable by the $A_n$-based test. As in the zero-mean case, the $A_n''$-based test outperforms Little’s test but is slightly outperformed by the $A_n$-based one. However, for the alternatives described in Remark \ref{remark_alg_mis}, the $A_n''$-based test attains the highest power, particularly at higher missingness rates; a representative example is shown in Figure \ref{fig:2X3Y_norm_mean1_Sigma05I05J_MAR_1_to_9_masked_r08}.

\subsection{Behavior as dimensionality increases}

Another important remark to make is that the standard implementation of Little’s MCAR test in the \texttt{R} package \texttt{naniar} can handle no more than 30 variables. To the best of our knowledge, the most capable implementation addressing this limitation is found in the (now deprecated) \texttt{Baylor}\-\texttt{EdPsych} package, which can process up to 50 variables. It is the structural flaw of any such test that splits the cases into groups with the same missingness pattern. For a fixed sample size, the more variables there are, the more patterns there will be, so the sample size for each possible eventually becomes equal to 1, and the test is rendered useless. In contrast, the novel test introduced here has no such constraints, neither theoretical nor practical. 

To examine the type I error and power behavior, we conduct a series of simulations. Following our previous notation, we generate 2X3Y, 5X5Y, and 10X10Y datasets and compare the performance of the tests as dimensionality increases.  The results for the missingness mechanisms described in the following paragraph are presented in the main text. Many others are placed in the Supplementary material.

MCAR data are generated in a standard way. MAR 1 to 9 missingness is constructed as follows. For 2X3Y data, missingness is generated using Algorithm \ref{algorithm_missingness} with $r = 0.5$. For the 5X5Y case, the algorithm was modified so that the variables $X^{(1)}$, $X^{(2)}$, and $X^{(3)}$ govern the missingness in $Y^{(1)}$, $Y^{(2)}$, and $Y^{(3)}$, respectively. Subsequently, the missingness in $Y^{(4)}$ and $Y^{(5)}$ is generated to be correlated with that in $Y^{(2)}$ and $Y^{(3)}$, following the same procedure as in the original algorithm. For the 10X10Y case, the same procedure was applied, with variables $X^{(1)}$ through $X^{(5)}$ governing the missingness in $Y^{(1)}$ through $Y^{(5)}$. The variables $Y^{(6)}$ through $Y^{(10)}$ were then made incomplete, with their response indicators correlated to those of $Y^{(1)}$ through $Y^{(5)}$. This choice was made so that the constucted alternatives are $A_n$-detectable, so we are able to evaluate the impact of using the more general $A_n''$-based test on power.

Simulation results for a normal distribution with all variable means equal to $1$, covariance matrix $0.5I + 0.5J$, and sample size $n = 100$, reveal that Little’s test suffers a substantial loss of both type I error control and power as dimensionality increases. This behavior is illustrated in Figure \ref{fig:largedim_sizes_n100} for type I error and in Figure \ref{fig:largedim powers_n100} for power. As dimensionality grows, the empirical type I error of Little’s test becomes much smaller than the nominal level, which negatively impacts its power.

In contrast, the $A_n''$-based test remains considerably more stable, with type I error showing noticeable inflation only for 20-dimensional data with high missingness rates. The $A_n$-based test is the most stable overall, with empirical type I error nearly equal to the nominal level. However, the improved $A_n''$-based test is preferable when no specific form of dependence among the response indicators can be assumed.

As shown in Figure \ref{fig:largedim_sizes_n300}, increasing the sample size to $n = 300$ allows the $A_n''$-based test to regain proper type I error control for 20-dimensional data, with empirical values remaining very close to the nominal level. In contrast, Little’s test continues to exhibit the same issues observed for smaller samples: its empirical type I error remains well below the nominal value, resulting in a marked reduction in power relative to the other two tests. 

 As shown in Figure \ref{fig:largedim powers_n100}, the generalized test exhibits lower power than the original $A_n$-based test, even in cases where type I error inflation occurs. The  power performance for $n=300$ is presented in Figure \ref{fig:largedim powers_n300}. The $A_n$-based test maintains consistent performance, while the $A_n''$-based test, though stable, has slightly lower power than the $A_n$-based one, similarly to the behavior observed in the smaller-sample case. Little's test is considerably behind both.

\begin{figure}
	\centering
	\includegraphics[width=\textwidth]{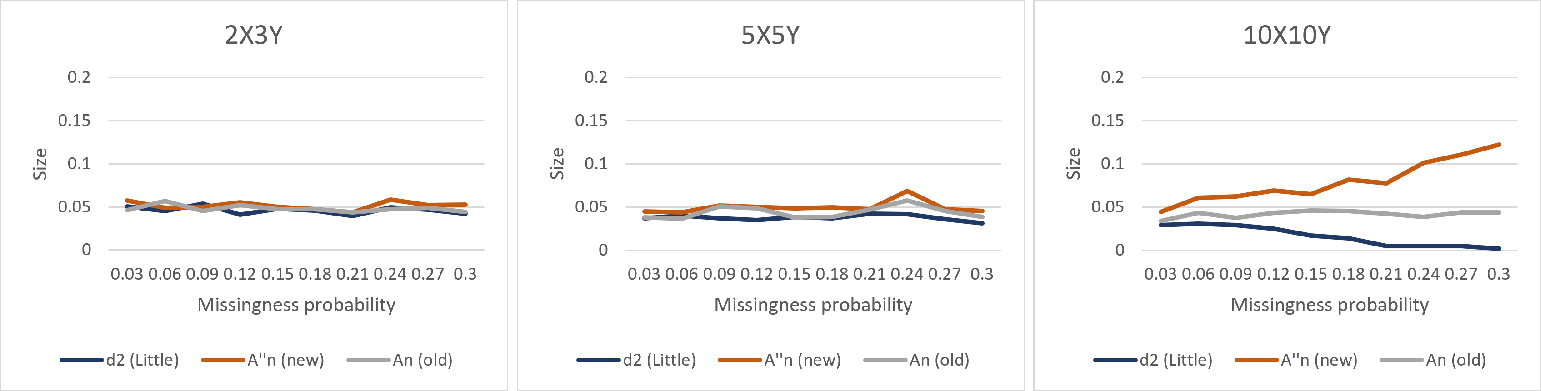}
	\caption{Empirical type I errors as dimension increases, normal distribution with mean $(1,\dots ,1)$ and covariance matrix $0.5I+0.5J$, MCAR with correlated response indicators, ${\rho} = 0.5$, $n = 100$}
	\label{fig:largedim_sizes_n100}
\end{figure}

\begin{figure}
	\centering
	\includegraphics[width=\textwidth]{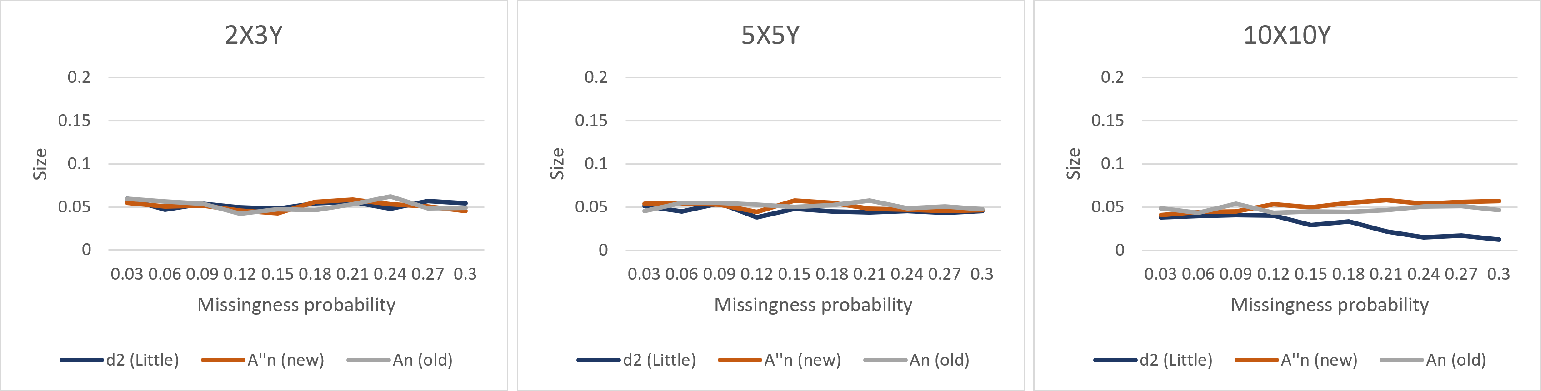}
	\caption{Empirical type I errors as dimension increases,  normal distribution with mean $(1,\dots ,1)$ and covariance matrix $0.5I+0.5J$, MCAR with correlated response indicators, ${\rho} = 0.5$, $n = 300$}
	\label{fig:largedim_sizes_n300}
\end{figure}

\begin{figure}
	\centering
	\includegraphics[width=\textwidth]{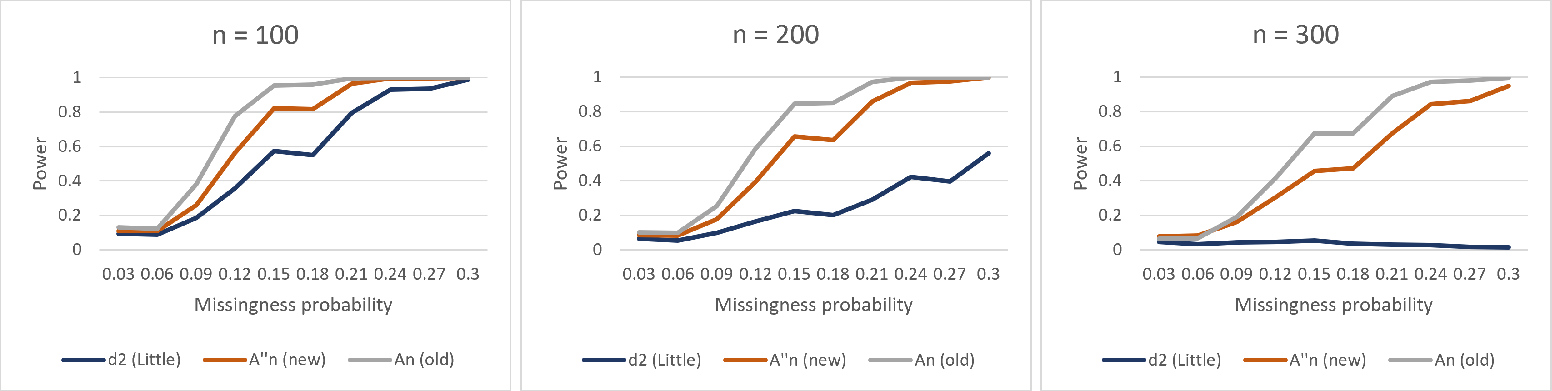}
	\caption{Empirical test powers as dimension increases, normal distribution with mean $(1,\dots ,1)$ and covariance matrix $0.5I+0.5J$, MAR 1 to 9 with correlated response indicators, ${\rho} = 0.5$, $n = 100$}
	\label{fig:largedim powers_n100}
\end{figure}

\begin{figure}
	\centering
	\includegraphics[width=\textwidth]{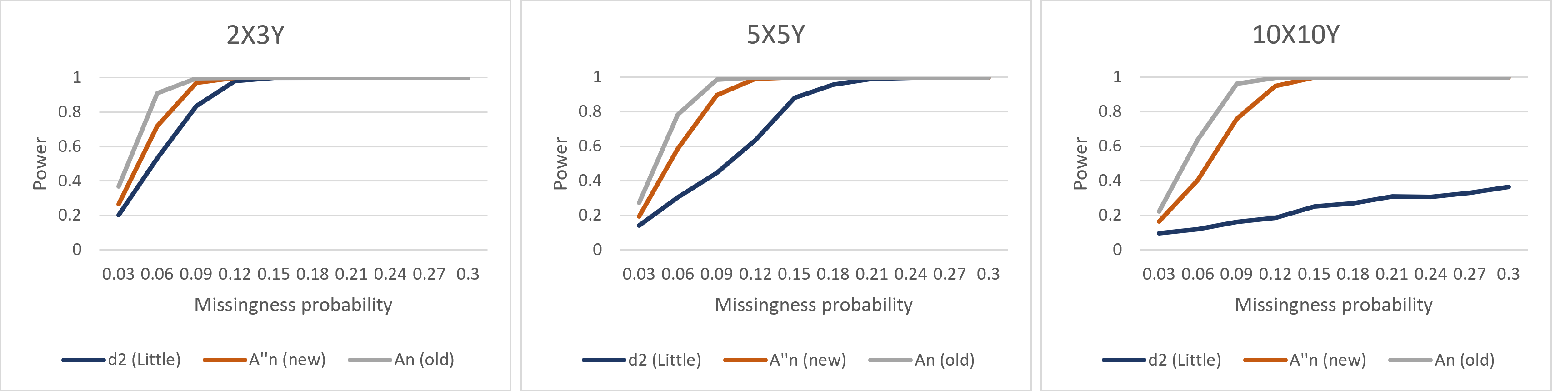}
	\caption{Empirical test powers as dimension increases, normal distribution with mean $(1,\dots ,1)$ and covariance matrix $0.5I+0.5J$, MAR 1 to 9 with correlated response indicators, ${\rho} = 0.5$, $n = 300$}
	\label{fig:largedim powers_n300}
\end{figure}


\section{Conclusion and outlook}\label{sec:conclusion}

The paper introduced an improved statistical test for assessing the MCAR assumption. Across all of the examined scenarios, particularly those that are more likely to arise in practice, such as cases with moderate missingness rates and large number of variables, the proposed test consistently outperformed Little’s MCAR test. In these settings, it demonstrated superior control of the type I error rate, higher statistical power, and greater robustness to {violations of the assumption of finite fourth moments}, provided that both tests performed satisfactorily. In situations involving infinite fourth moments combined with the alternatives that are more difficult to detect, both the novel test and Little’s test exhibited unexpected behavior. Specifically, their power declined as the missingness rate increased, which {was unexpected}. 

However, in contrast to Little’s test, the novel test did not exhibit a substantial loss of power as dimensionality increased, indicating that it performs more reliably in high-dimensional settings. This stability suggests that the novel test may be better suited for modern applications involving large number of variables relative to the sample size, where traditional methods often struggle. A natural direction for future research would be to investigate the asymptotic properties of the proposed test, as well as Little’s test, as the dimension grows. In particular, it would be of interest to derive their asymptotic distributions when the dimension tends to infinity: either at the same rate as, faster than, or slower than the sample size.

With regard to Remark \ref{remark_transform}, another possible direction for further improvement would be to replace the covariance, which is used in the current formulation as a measure of linear dependence, with an alternative discrepancy measure that {either characterizes dependence or is more closely related to it}. {This is a potential goal for future research.}

\cite{heinze2024phases} proposed four phases of methodological research that, although developed {primarily} for the biostatistical framework, are broadly relevant across statistics. Their brief overview is as follows. Phase I involves the theoretical development of a new method. Phase II focuses on empirical evaluation in a narrow setting. Phase III includes validation across diverse scenarios and the creation of {user-friendly} software implementations of the proposed methods. Phase IV aims for comprehensive understanding of the method, including knowing when it is preferred or not, identifying common pitfalls, and developing practical diagnostics of whether the assumptions of the method are met.

Our work {currently} falls between Phase II and Phase III.  Two main goals for future work in this context would be to broaden the simulation study such that more variability is added to sample size, data distribution, dimension, as well as the underlying distribution. A very important class of distributions to study the test's performance on are discrete distributions, especially those of nominal type (yes/no, male/female etc.) It would be interesting to see whether the usage of different encodings has an impact on the power performance of the test. It seems to us that those types of simulations fall outside of the scope of this paper at this moment. 

Comparing the proposed test to other tests of MCAR, that are not as widely used as Little's test, is also an interesting topic. Currently, the literature is lacking a paper with systematic comparison of existing MCAR tests in terms of empirical type I error and power.

Finally, phase IV of methodological research involves gaining a deeper understanding of the method through practical use. We hope and expect that as our method is {adopted} more widely, both its strengths and limitations will become clearer.

\section*{Supplementary material}

\texttt{R} implementation of function that return the $p$-values of both $A_n$- and $A_n''$-based tests, the code used for running the simulations, as well as additional simulation results, can be found on the author's \texttt{GitHub} profile: \url{https://github.com/danijel-g-aleksic}, or on request.

\section*{Acknowledgements}

The author would like to express sincere gratitude to the anonymous referee of his original paper \citep{aleksic2024novel}, who insisted there must be a way to include the partially observed data into the framework of covariance-based test. Their criticisms initiated the thought process that lead to the creation of this paper. 

Next, the author is grateful to the (other) anonymous referee of the same paper, who pointed out the possibility shown in Remark \ref{remark_transform}.

The author would also like to thank professor Bojana Milošević, PhD, from University of Belgrade -- Faculty of Mathematics, for a series of useful remarks that improved the quality and the structure of this paper. 

\section*{Disclosure statement}

No potential conflict of interest was reported by the author.

\section*{Funding}

The work of D.G. Aleksić is supported by the Ministry of Science, Technological Development and Innovations of the Republic of Serbia (the contract 451-03-137/2025-03/200151).

\bibliographystyle{abbrvnat}
\bibliography{literatura}

@article{Hoeffding1948,
author = {Wassily Hoeffding},
title = {{A Class of Statistics with Asymptotically Normal Distribution}},
volume = {19},
journal = {The Annals of Mathematical Statistics},
number = {3},
publisher = {Institute of Mathematical Statistics},
pages = {293 -- 325},
year = {1948},
doi = {10.1214/aoms/1177730196}
}

@book{serfling1980approximation,
  title={{Approximation Theorems of Mathematical Statistics}},
  author={Serfling, Robert J},
  year={1980},
  publisher={John Wiley \& Sons},
address={{New York}},
doi={10.1002/9780470316481}
}

@book{LittleRubin1987,
  author = {Roderick J.A. Little and Donald B. Rubin},
  year = {1987},
  title = {{Statistical Analysis with Missing Data, First Edition}},
  publisher = {John Wiley \& Sons},
    address={New York}
}

@article{Little1988,
author = {Roderick J. A. Little},
title = {{A Test of Missing Completely at Random for Multivariate Data with Missing Values}},
journal = {Journal of the American Statistical Association},
volume = {83},
number = {404},
pages = {1198-1202},
year = {1988},
publisher = {Taylor & Francis},
DOI = {10.1080/01621459.1988.10478722},
}

@article{santos2019generating,
  title={{Generating synthetic missing data: A review by missing mechanism}},
  author={Santos, Miriam Seoane and Pereira, Ricardo Cardoso and Costa, Adriana Fonseca and Soares, Jastin Pompeu and Santos, Jo{\~a}o and Abreu, Pedro Henriques},
  journal={IEEE Access},
  volume={7},
  pages={11651--11667},
  year={2019},
  publisher={IEEE},
doi={10.1109/ACCESS.2019.2891360}
}

@article{aleksic2023etAl,
  title={{Non-degenerate U-statistics for data missing completely at random with application to testing independence}},
  author={Aleksi{\'c}, Danijel and Cupari{\'c}, Marija and Milo{\v{s}}evi{\'c}, Bojana},
  journal={Stat},
  volume={12},
  number={1},
  pages={e634},
  year={2023},
  publisher={Wiley Online Library},
    doi = {10.1002/sta4.634}
}

@Manual{missMethods,
    title = {{missMethods: Methods for Missing Data}},
    author = {Tobias Rockel},
    year = {2023},
    note = {R package version 0.4.0},
    url = {https://CRAN.R-project.org/package=missMethods},
  }

@article{li2015nonparametric,
  title={{A nonparametric test of missing completely at random for incomplete multivariate data}},
  author={Li, Jun and Yu, Yao},
  journal={Psychometrika},
  volume={80},
  pages={707--726},
  year={2015},
  publisher={Springer},
doi={10.1007/s11336-014-9410-4}
}

@article{fuchs1982maximum,
  title={{Maximum likelihood estimation and model selection in contingency tables with missing data}},
  author={Fuchs, Camil},
  journal={Journal of the American Statistical Association},
  volume={77},
  number={378},
  pages={270--278},
  year={1982},
  publisher={Taylor \& Francis},
doi={10.2307/2287230}
}

@article{jamshidian2010tests,
  title={{Tests of homoscedasticity, normality, and missing completely at random for incomplete multivariate data}},
  author={Jamshidian, Mortaza and Jalal, Siavash},
  journal={Psychometrika},
  volume={75},
  number={4},
  pages={649--674},
  year={2010},
  publisher={Springer},
doi={10.1007/s11336-010-9175-3}
}

@article{kim2002tests,
  title={{Tests of homogeneity of means and covariance matrices for multivariate incomplete data}},
  author={Kim, Kevin H and Bentler, Peter M},
  journal={Psychometrika},
  volume={67},
  pages={609--623},
  year={2002},
  publisher={Springer},
doi={10.1007/BF02295134}
}

@article{spohn2021pklm,
  title={{PKLM: A flexible MCAR test using Classification}},
  author={Spohn, Meta-Lina and N{\"a}f, Jeffrey and Michel, Loris and Meinshausen, Nicolai},
  journal={arXiv preprint arXiv:2109.10150},
  year={2021},
doi={10.48550/arXiv.2109.10150}
}

@article{berrett2023optimal,
  title={{Optimal nonparametric testing of Missing Completely At Random and its connections to compatibility}},
  author={Berrett, Thomas B and Samworth, Richard J},
  journal={The Annals of Statistics},
  volume={51},
  number={5},
  pages={2170--2193},
  year={2023},
  publisher={Institute of Mathematical Statistics},
doi={10.1214/23-AOS2326}
}

@article{bordino2025tests,
	title={Tests of missing completely at random based on sample covariance matrices},
	author={Bordino, Alberto and Berrett, Thomas B},
	journal={The Annals of Statistics},
	volume={53},
	number={5},
	pages={2204--2229},
	year={2025},
	publisher={Institute of Mathematical Statistics}
}

@article{aleksic2024novel,
  title={{A Novel Test of Missing Completely at Random: U-statistics-based Approach}},
  author={Aleksi{\'c}, Danijel},
  journal={Statistics},
  year={2024},
    volume={51},
  number={5},
  pages={2170--2193},
    doi = {10.1080/02331888.2024.2386361}
}

@article{diggle1989testing,
  title={{Testing for random dropouts in repeated measurement data}},
  author={Diggle, Peter J},
  journal={Biometrics},
  pages={1255--1258},
volume={45},
number={4},
  year={1989},
doi={10.2307/2531777}
}

@article{ridout1991testing,
  title={{Testing for random dropouts in repeated measurement data}},
  author={Ridout, Martin S and Diggle, Peter J},
  journal={Biometrics},
  pages={1617--1621},
volume={47},
number={4},
  year={1991},
doi={10.2307/2532413}
}

@article{park1993atest,
  title={{A test of the missing data mechanism for repeated categorical data}},
  author={Park, Taesung and Davis, Charles S},
  journal={Biometrics},
volume={49},
number={2},
  pages={631--638},
  year={1993},
doi={10.2307/2532576}
}

@article{park1993btest,
  title={{A test of the missing data mechanism for repeated measures data}},
  author={Park, Taesung and Lee, Seungyeoun and Woolson, Robert F},
  journal={Communications in Statistics-Theory and Methods},
  volume={22},
  number={10},
  pages={2813--2829},
  year={1993},
  publisher={Taylor \& Francis},
doi={10.1080/03610929308831187}
}

@article{park1997test,
  title={{A test of missing completely at random for longitudinal data with missing observations}},
  author={Park, Taesung and Lee, Seung-Yeoun},
  journal={{Statistics in Medicine}},
  volume={16},
  number={16},
  pages={1859--1871},
  year={1997},
  publisher={Wiley Online Library},
doi={10.1002/(SICI)1097-0258(19970830)16:16<1859::AID-SIM593>3.0.CO;2-3}
}

@article{listing1998tests,
  title={{Tests if dropouts are missed at random}},
  author={Listing, Joachim and Schlittgen, Rainer},
  journal={Biometrical Journal: Journal of Mathematical Methods in Biosciences},
  volume={40},
  number={8},
  pages={929--935},
  year={1998},
  publisher={Wiley Online Library},
doi={10.1002/(SICI)1521-4036(199812)40:8<929::AID-BIMJ929>3.0.CO;2-X}
}

@article{chen1999test,
  title={{A test of missing completely at random for generalised estimating equations with missing data}},
  author={Chen, Hua Yun and Little, Roderick},
  journal={Biometrika},
  volume={86},
  number={1},
  pages={1--13},
  year={1999},
  publisher={Oxford University Press},
doi={10.1093/biomet/86.1.1}
}

@article{qu2002testing,
  title={{Testing ignorable missingness in estimating equation approaches for longitudinal data}},
  author={Qu, Annie and Song, Peter X-K},
  journal={Biometrika},
  volume={89},
  number={4},
  pages={841--850},
  year={2002},
  publisher={Oxford University Press},
doi={10.1093/biomet/89.4.841}
}

@techreport{fairclough2003design,
  title={{Design and analysis of quality of life studies in clinical trials}},
  author={Fairclough, Diane Lynn},
  year={2002},
  institution={CRC press},
ISBN={1-58488-263-8}
}

@article{potthoff2006can,
  title={{Can one assess whether missing data are missing at random in medical studies?}},
  author={Potthoff, Richard F and Tudor, Gail E and Pieper, Karen S and Hasselblad, Vic},
  journal={{Statistical Methods in Medical Research}},
  volume={15},
  number={3},
  pages={213--234},
  year={2006},
  publisher={Sage Publications Sage CA: Thousand Oaks, CA},
doi={10.1191/0962280206sm448oa}
}

@article{jamshidian2007testing,
  title={{Testing equality of covariance matrices when data are incomplete}},
  author={Jamshidian, Mortaza and Schott, James R},
  journal={Computational Statistics \& Data Analysis},
  volume={51},
  number={9},
  pages={4227--4239},
  year={2007},
  publisher={Elsevier},
doi={10.1016/j.csda.2006.05.005}
}

@article{jamshidian2008postmodeling,
  title= {{Postmodeling sensitivity analysis to detect the effect of missing data mechanisms}},
  author={Jamshidian, Mortaza and Mata, Matthew},
  journal={Multivariate Behavioral Research},
  volume={43},
  number={3},
  pages={432--452},
  year={2008},
  publisher={Taylor \& Francis},
doi={10.1080/00273170802285792}
}

@article{fielding2009investigating,
  title={{Investigating the missing data mechanism in quality of life outcomes: a comparison of approaches}},
  author={Fielding, Shona and Fayers, Peter M and Ramsay, Craig R},
  journal={Health and Quality of Life Outcomes},
  volume={7},
  pages={1--10},
  year={2009},
  publisher={Springer},
doi={10.1186/1477-7525-7-57}
}

@article{jamshidian2013data,
  title={{Data-driven sensitivity analysis to detect missing data mechanism with applications to structural equation modelling}},
  author={Jamshidian, Mortaza and Yuan, Ke-Hai},
  journal={Journal of Statistical Computation and Simulation},
  volume={83},
  number={7},
  pages={1344--1362},
  year={2013},
  publisher={Taylor \& Francis},
doi={10.1080/00949655.2012.660486}
}

@book{lin2013probability,
  title=  {{A probability based framework for testing the missing data mechanism {(PhD Thesis)}}},
  author={Lin, Johnny Cheng-Han},
  year={2013},
  publisher={University of California},
address={Los Angeles},
URL = {https://escholarship.org/content/qt4c51m4bm/qt4c51m4bm_noSplash_7fb298286f062237f477d811c9d95f11.pdf}
}

@article{jamshidian2014examining,
  title={{Examining missing data mechanisms via homogeneity of parameters, homogeneity of distributions, and multivariate normality}},
  author={Jamshidian, Mortaza and Yuan, Ke-Hai},
  journal={Wiley Interdisciplinary Reviews: Computational Statistics},
  volume={6},
  number={1},
  pages={56--73},
  year={2014},
  publisher={Wiley Online Library},
doi={10.1002/wics.1287}
}

@article{yuan2018missing,
  title={{Missing data mechanisms and homogeneity of means and variances--covariances}},
  author={Yuan, Ke-Hai and Jamshidian, Mortaza and Kano, Yutaka},
  journal={Psychometrika},
  volume={83},
  pages={425--442},
  year={2018},
  publisher={Springer},
doi={10.1007/s11336-018-9609-x}
}

@article{zhang2019unified,
  title={{A unified empirical likelihood approach for testing MCAR and subsequent estimation}},
  author={Zhang, Shixiao and Han, Peisong and Wu, Changbao},
  journal={Scandinavian Journal of Statistics},
  volume={46},
  number={1},
  pages={272--288},
  year={2019},
  publisher={Wiley Online Library},
doi={10.1111/sjos.12351}
}

@article{bojinov2020diagnosing,
  title={{Diagnosing missing always at random in multivariate data}},
  author={Bojinov, Iavor I and Pillai, Natesh S and Rubin, Donald B},
  journal={Biometrika},
  volume={107},
  number={1},
  pages={246--253},
  year={2020},
  publisher={Oxford University Press},
doi={10.1093/biomet/asz061}
}

@article{rouzinov2022regression,
  title={{Regression-based approach to test missing data mechanisms}},
  author={Rouzinov, Serguei and Berchtold, Andr{\'e}},
  journal={Data},
  volume={7},
  number={2},
  pages={1--16},
  year={2022},
  publisher={MDPI},
doi={10.3390/data7020016}
}

@article{chassan2023test,
  title={{How to test the missing data mechanism in a hidden Markov model}},
  author={Chassan, Malika and Concordet, Didier},
  journal={Computational Statistics \& Data Analysis},
  volume={182},
  pages={107--723},
  year={2023},
  publisher={Elsevier},
doi={10.1016/j.csda.2023.107723}
}

@article{aleksic2024impute,
  title={{To impute or not? Testing multivariate normality on incomplete dataset: revisiting the {BHEP} test}},
  author={Aleksi{\'c}, Danijel G and Milo{\v{s}}evi{\'c}, Bojana},
  journal={Journal of Applied Statistics},
 volume={52},
  number={9},
  pages={1742--1759},
  year={2025},
  publisher={Taylor \& Francis},
doi={10.1080/02664763.2024.2438798}
}

@article{cuparic2024ipcw,
  title={{IPCW approach for testing independence}},
  author={Cupari{\'c}, Marija and Milo{\v{s}}evi{\'c}, Bojana},
  journal={Journal of Nonparametric Statistics},
  volume={36},
  number={1},
  pages={118--145},
  year={2024},
  publisher={Taylor \& Francis},
doi={10.1080/10485252.2023.2185749}
}

@article{heinze2024phases,
  title={{Phases of methodological research in biostatistics—building the evidence base for new methods}},
  author={Heinze, Georg and Boulesteix, Anne-Laure and Kammer, Michael and Morris, Tim P and White, Ian R and Simulation Panel of the STRATOS Initiative},
  journal={Biometrical Journal},
  volume={66},
  number={1},
  pages={2200222},
  year={2024},
  publisher={Wiley Online Library},
doi={10.1002/bimj.202200222}
}

@article{ofner2025testing,
  title={{Testing the Missing Completely at Random Assumption for Functional Data}},
  author={Ofner, Maximilian and H{\"o}rmann, Siegfried and Kraus, David and Liebl, Dominik},
  journal={arXiv preprint arXiv:2505.08721},
  year={2025},
doi={10.48550/arXiv.2505.08721}
}

@article{sklar1959fonctions,
  title={{Fonctions de r{\'e}partition {\`a} n dimensions et leurs marges}},
  author={Sklar, A},
  journal={Annales de l'ISUP},
  volume={8},
  number={3},
  pages={229--231},
  year={1959}
}

@online{github,
  author = {Aleksi{\'c}, Danijel G.},
  title = {{\texttt{GitHub}} page:  \url{https://github.com/danijel-g-aleksic}},
  year = 2025
}

@article{kochar1990distribution,
  title={Distribution-free tests based on sub-sample extrema for testing against positive dependence},
  author={Kochar, Subhash C and Gupta, RP},
  journal={Australian Journal of Statistics},
  volume={32},
  number={1},
  pages={45--51},
  year={1990},
  publisher={Wiley Online Library}
}

@article{fischer2012constructing,
  title={{Constructing and generalizing given multivariate copulas: A unifying approach}},
  author={Fischer, Matthias and K{\"o}ck, Christian},
  journal={Statistics},
  volume={46},
  number={1},
  pages={1--12},
  year={2012},
  publisher={Taylor \& Francis}
}

@phdthesis{aleksic2026PhD,
  author      = {Aleksić, Danijel G.},
  title       = {U- and {V}- statistics for incomplete data and their application to model specification testing},
  school      = {University of Belgrade, Faculty of Mathematics},
  address     = {Belgrade, Serbia},
  year        = {2026}
}

@book{van2000asymptotic,
  title={Asymptotic statistics},
  author={Van der Vaart, Aad W},
  volume={3},
  year={2000},
  publisher={Cambridge university press}
}

\appendix

\section{Proofs}\label{Appendix_proofs}

\subsection{Proof of Lemma \ref{lemma_YY}}\label{lemma_YY_proof}

Let $u,v,r,s$ be fixed. Begin by noting that
\begin{align*}
    \hat{T}_{n, Y}^{(u,v)} = \frac{1}{\binom{n}{1}} \sum_{i=1}^n Y_i^{(u)}R_i^{(u)}R_i^{(v)} - \frac{1}{\binom{n}{2}}\underset{1 \leq i < j \leq n}{\sum \sum } \frac{1}{2} \bigg( Y_i^{(u)}R_i^{(u)}R_j^{(v)} + Y_j^{(u)}R_j^{(u)}R_i^{(v)}\bigg) =: M_n - N_n
\end{align*}
and, similarly
\begin{align*}
    \hat{T}_{n, Y}^{(r,s)} = \frac{1}{\binom{n}{1}} \sum_{i=1}^n Y_i^{(r)}R_i^{(r)}R_i^{(s)} - \frac{1}{\binom{n}{2}}\underset{1 \leq i < j \leq n}{\sum \sum } \frac{1}{2} \bigg( Y_i^{(r)}R_i^{(r)}R_j^{(s)} + Y_j^{(r)}R_j^{(r)}R_i^{(s)}\bigg) =: Q_n - S_n.
\end{align*}
Relying on Theorem 7.1 from \cite{Hoeffding1948}, we have that
\begin{align*}
    \lim_{n \to \infty} n \Cov \left( M_n, Q_n \right) &= 1\cdot 1 \cdot\Cov \left( Y_1^{(u)} R_1^{(u)} R_1^{(v)},\, Y_1^{(r)} R_1^{(r)} R_1^{(s)} \right) \\
    &= \E \left( Y^{(u)}Y^{(r)} \right) \E \left( R^{(u)}R^{(v)} R^{(r)} R^{(s)} \right) - \E \left( Y^{(u)} \right) \E \left( Y^{(r)}\right) \E \left(R^{(u)}R^{(v)} \right) \E \left( R^{(r)} R^{(s)} \right),
\end{align*}
and, similarly,
\begin{align*}
    \lim_{n \to \infty}& n \Cov \left( M_n, S_n\right) \\
    &= 1 \cdot 2 \cdot \Cov \left( Y_1^{(u)} R_1^{(u)}R_1^{(v)}, \frac12 \left(  Y_1^{(r)} R_1^{(r)}R_2^{(s)} + Y_2^{(r)}R_2^{(r)} R_1^{(s)}  \right) \right) \\
    &= \Cov \left( Y_1^{(u)}R_1^{(u)} R_1^{(v)},\, Y_1^{(r)}R_1^{(r)} R_2^{(s)} \right) + \Cov \left( Y_1^{(u)}R_1^{(u)} R_1^{(v)}, Y_2^{(r)}R_2^{(r)}R_1^{(s)} \right) \\
    &= \E \left( Y^{(u)}Y^{(r)} \right) \E \left( R^{(u)}R^{(v)}R^{(r)} \right) \E \left( R^{(s)} \right) - \E \left( Y^{(u)} \right) \E \left( Y^{(r)} \right) \E \left( R^{(u)} R^{(v)} \right) \E \left( R^{(r)} \right) \E \left( R^{(s)} \right) \\
    & \qquad + \E \left( Y^{(u)} \right) \E \left( Y^{(r)} \right)  \E \left( R^{(u)} R^{(v)} R^{(s)} \right) \E \left( R^{(r)} \right) - \E \left( Y^{(u)} \right) \E \left( Y^{(r)} \right) \E \left( R^{(u)} R^{(v)} \right) \E \left( R^{(r)} \right) \E \left( R^{(s)} \right) \\
    &= \E \left( Y^{(u)}Y^{(r)} \right) \E \left( R^{(u)}R^{(v)}R^{(r)} \right) \E \left( R^{(s)} \right) -2 \E \left( Y^{(u)} \right) \E \left( Y^{(r)} \right) \E \left( R^{(u)} R^{(v)} \right) \E \left( R^{(r)} \right) \E \left( R^{(s)} \right) \\
    & \qquad + \E \left( Y^{(u)} \right) \E \left( Y^{(r)} \right)  \E \left( R^{(u)} R^{(v)} R^{(s)} \right) \E \left( R^{(r)} \right).
\end{align*}
Analogously, we obtain
\begin{align*}
    \lim_{n \to \infty}& n \Cov \left( N_n, Q_n\right) \\
    &= 2 \cdot 1 \cdot \Cov \left( \frac12 \left( Y_1^{(u)}R_1^{(u)}R_2^{(v)} + Y_2^{(u)}R_2^{(u)}R_1^{(v)} \right),\, Y_1^{(r)}R_1^{(r)}R_1^{(s)}  \right)\\
    &= \Cov \left( Y_1^{(u)}R_1^{(u)}R_2^{(v)}, Y_1^{(r)} R_1^{(r)}R_1^{(s)} \right) + \Cov \left( Y_2^{(u)}R_2^{(u)} R_1^{(v)}, Y_1^{(r)} R_1^{(r)}R_1^{(s)}  \right) \E \left(R^{(v)} \right)   \\
    &= \E \left( Y^{(u)}Y^{(r)} \right) \E \left( R^{(u)}R^{(r)}R^{(s)} \right) + \E \left( Y^{(u)} \right) \E \left( Y^{(r)} \right) \E \left(R^{(v)}R^{(r)}R^{(s)} \right)\E \left(R^{(u)} \right) \\
    & \qquad -2 \E \left( Y^{(u)} \right) \E \left( Y^{(r)}\right) \E \left(R^{(u)}R^{(v)} \right) \E \left( R^{(r)}R^{(s)}\right)
\end{align*}
and
\begin{align*}
    \lim_{n \to \infty}& n \Cov \left( N_n, S_n\right) \\
    &= 2 \cdot 2 \cdot \Cov \left( \frac12 \left( Y_1^{(u)}R_1^{(u)}R_2^{(v)} +   Y_2^{(u)}R_2^{(u)}R_1^{(v)} \right), \, \frac12 \left( Y_1^{(r)}R_1^{(r)}R_3^{(s)} +   Y_3^{(r)}R_3^{(r)}R_1^{(s)} \right) \right) \\
    & = \Cov \left( Y_1^{(u)}R_1^{(u)}R_2^{(v)}, \, Y_1^{(r)}R_1^{(r)}R_3^{(s)} \right) + \Cov \left( Y_1^{(u)}R_1^{(u)}R_2^{(v)}, \, Y_3^{(r)}R_3^{(r)}R_1^{(s)} \right) \\
    & \qquad + \Cov \left( Y_2^{(u)}R_2^{(u)}R_1^{(v)}, \, Y_1^{(r)}R_1^{(r)}R_3^{(s)} \right) + \Cov \left( Y_2^{(u)}R_2^{(u)}R_1^{(v)}, \, Y_3^{(r)}R_3^{(r)}R_1^{(s)} \right) \\
    &= \E \left( Y^{(u)}Y^{(r)}\right) \E \left( R^{(u)}R^{(r)} \right) \E \left( R^{(v)}  \right) \E \left(  R^{(s)}\right) + \E \left( Y^{(u)} \right) \E \left( Y^{(r)} \right) \E \left( R^{(u)}R^{(s)} \right) \E \left(  R^{(v)}\right) \E \left(  R^{(r)}\right) \\
    & \qquad + \E \left( Y^{(u)} \right) \E \left(Y^{(r)}  \right) \E \left(  R^{(v)}R^{(r)}\right) \E \left( R^{(u)} \right) \E \left( R^{(s)} \right) + \E \left(Y^{(u)}  \right) \E \left( Y^{(r)} \right) \E \left( R^{(v)}R^{(s)}  \right) \E \left( R^{(u)} \right) \E \left( R^{(r)} \right) \\
    & \qquad - 4 \E \left( Y^{(u)} \right) \E \left( Y^{(r)} \right) \E \left(  R^{(u)}\right)\E \left( R^{(v)} \right) \E \left( R^{(r)} \right) \E \left( R^{(s)} \right).
\end{align*}
By noting that
\begin{align*}
    \lim_{n \to \infty} n \Cov \left( \hat{T}_{n,Y}^{(u,v)}, \, \hat{T}_{n,Y}^{(r,s)} \right) &= \lim_{n \to \infty} n \Cov \left( M_n, Q_n \right) - \lim_{n \to \infty} n \Cov \left( M_n, S_n \right) \\
    & \qquad \qquad - \lim_{n \to \infty} n \Cov \left( N_n, Q_n \right) + \lim_{n \to \infty} n \Cov \left( N_n, S_n \right)
\end{align*}
and combining the derived expressions, we obtain the statement of the Lemma. This concludes the proof.

\subsection{Proof of Corollary \ref{corollary_XY}}\label{proof_corollary_XY}

The results follows from Lemma \ref{corollary_XY} by setting $Y^{(u)} = X^{(u)}$ and noting that $R^{(u)} \equiv 1$. 

\subsection{Proof of Corollary \ref{corollary_XX}}\label{proof_corollary_XX}

Follows directly from Corollary \ref{corollary_XY} by setting $Y^{(r)} = X^{(r)}$ and noting that $R^{(r)} \equiv 1$.

\subsection{Proof of Theorem \ref{theorem_main}}\label{proof_main}

The results follows directly from the equations \eqref{cov_YY_omitted}, \eqref{cov_XY}, and \eqref{cov_XX}, and the fact that the $\hat{\Lambda}$ is a consistent estimator under the finite fourth moments assumption.

\subsection{Proof of Theorem \ref{theorem_plugin}}\label{proof_theorem_plugin}

Let $\tilde{Y}_i^{(u)} = Y_i^{(u)} - \mu^{(u)}$ define the truly centered values. We can write the sample-centered statistic explicitly as a function of the plug-in vector $\boldsymbol{\hat{\mu}}_Y$:
\begin{align*}
    \hat{T}_{n, Y}^{(u, v)}(\boldsymbol{\hat{\mu}}_Y) &= \frac{1}{n} \sum_{i=1}^n \left(Y_i^{(u)} - \hat{\mu}_Y^{(u)}\right) R_i^{(u)} R_i^{(v)} - \frac{1}{n(n-1)} \sum_{i=1}^n \sum_{\substack{j=1 \\ j \neq i}}^n \left(Y_i^{(u)} - \hat{\mu}_Y^{(u)}\right) R_i^{(u)} R_j^{(v)} \\
    &= \hat{T}_{n, Y}^{(u, v)}(\boldsymbol{\mu}_Y) - \left(\hat{\mu}_Y^{(u)} - \mu^{(u)}\right) \left[ \frac{1}{n}\sum_{i=1}^n R_i^{(u)}R_i^{(v)} - \frac{1}{n(n-1)}\sum_{i=1}^n \sum_{\substack{j=1 \\ j \neq i}}^n R_i^{(u)}R_j^{(v)} \right].
\end{align*}
Let the bracketed expression be denoted as $G_n^{(u,v)}$. By the weak law of large numbers for $U$-statistics, $G_n^{(u,v)}$ converges in probability to its expected value under MCAR:
\begin{align*}
    G_n^{(u,v)} \overset{P}{\to} \E\left(R^{(u)}R^{(v)}\right) - \E\left(R^{(u)}\right)\E\left(R^{(v)}\right) = \Cov\left(R^{(u)}, R^{(v)}\right).
\end{align*}
Thus, applying Slutsky's theorem, we can express the scaling behavior as:
\begin{align}\label{slutsky_linearization}
    \sqrt{n}\hat{T}_{n, Y}^{(u, v)}(\boldsymbol{\hat{\mu}}_Y) = \sqrt{n}\hat{T}_{n, Y}^{(u, v)}(\boldsymbol{\mu}_Y) - \Cov\left(R^{(u)}, R^{(v)}\right) \sqrt{n}\left(\hat{\mu}_Y^{(u)} - \mu^{(u)}\right) + o_P(1).
\end{align}

Under the MCAR null hypothesis, the standard first-order Taylor expansion for the available-case mean yields the influence representation:
\begin{align}\label{mean_influence}
    \sqrt{n}\left(\hat{\mu}_Y^{(u)} - \mu^{(u)}\right) = \frac{1}{\sqrt{n} \E\left(R^{(u)}\right)} \sum_{i=1}^n \tilde{Y}_i^{(u)} R_i^{(u)} + o_P(1).
\end{align}

To compute the first-order Hoeffding projection of the true-mean statistic $\hat{T}_{n, Y}^{(u, v)}(\boldsymbol{\mu}_Y)$, we first express it as a single, unified $U$-statistic of degree 2. By representing the single-sum component over a symmetric double-sum index space, the statistic can be rewritten as:
\begin{align*}
    \hat{T}_{n, Y}^{(u, v)}(\boldsymbol{\mu}_Y) = \frac{1}{\binom{n}{2}} \sum_{1 \leq i < j \leq n} h\left(O_i, O_j\right),
\end{align*}
where $O_i = (\tilde{Y}_i^{(u)}, R_i^{(u)}, R_i^{(v)})$ represents the observed vector for the $i$-th case, and the symmetric kernel $h(O_i, O_j)$ simplifies to the elegant product form:
\begin{align*}
    h\left(O_i, O_j\right) = \frac{1}{2} \left( \tilde{Y}_i^{(u)}R_i^{(u)} - \tilde{Y}_j^{(u)}R_j^{(u)} \right) \left( R_i^{(v)} - R_j^{(v)} \right).
\end{align*}

Following standard $U$-statistic theory \citep[][Chapter 5]{serfling1980approximation}, the first-order projection of a degree-2 $U$-statistic is defined as $\frac{2}{n}\sum_{i=1}^n h_1(O_i)$, where $h_1(O_i) = \E[h(O_i, O_j) \mid O_i]$. Expanding the product kernel inside this conditional expectation gives:
\begin{align*}
    h_1(O_i) &= \frac{1}{2} \E \left[ \tilde{Y}_i^{(u)}R_i^{(u)}R_i^{(v)} - \tilde{Y}_i^{(u)}R_i^{(u)}R_j^{(v)} - \tilde{Y}_j^{(u)}R_j^{(u)}R_i^{(v)} + \tilde{Y}_j^{(u)}R_j^{(u)}R_j^{(v)} \;\middle|\; O_i \right].
\end{align*}
Under the MCAR null hypothesis, the data values and the response indicators are independent. Since the variables are centered ($\E(\tilde{Y}^{(u)}) = 0$), any terms containing an unconditioned $\tilde{Y}_j^{(u)}$ vanish completely because $\E[\tilde{Y}_j^{(u)}R_j^{(u)}] = 0$ and $\E[\tilde{Y}_j^{(u)}R_j^{(u)}R_j^{(v)}] = 0$. This leaves:
\begin{align*}
    h_1(O_i) = \frac{1}{2} \tilde{Y}_i^{(u)} R_i^{(u)} \left( R_i^{(v)} - \E\left(R^{(v)}\right) \right).
\end{align*}

Since a non-degenerate $U$-statistic is asymptotically equivalent to its first-order projection under $\sqrt{n}$-scaling, we combine the terms and obtain the linear approximation:
\begin{align}\label{hoeffding_projection_true}
    \sqrt{n}\hat{T}_{n, Y}^{(u, v)}(\boldsymbol{\mu}_Y) = \frac{1}{\sqrt{n}} \sum_{i=1}^n \tilde{Y}_i^{(u)} R_i^{(u)} \left( R_i^{(v)} - \E\left(R^{(v)}\right) \right) + o_P(1).
\end{align}

Substituting \eqref{mean_influence} and \eqref{hoeffding_projection_true} back into the linearized relationship \eqref{slutsky_linearization} allows us to merge the influence components:
\begin{align*}
    \sqrt{n}\hat{T}_{n, Y}^{(u, v)}(\boldsymbol{\hat{\mu}}_Y) &= \frac{1}{\sqrt{n}} \sum_{i=1}^n \left[ \tilde{Y}_i^{(u)} R_i^{(u)} \left( R_i^{(v)} - \E\left(R^{(v)}\right) \right) - \frac{\Cov\left(R^{(u)}, R^{(v)}\right)}{\E\left(R^{(u)}\right)} \tilde{Y}_i^{(u)} R_i^{(u)} \right] + o_P(1) \\
    &= \frac{1}{\sqrt{n}} \sum_{i=1}^n \tilde{Y}_i^{(u)} R_i^{(u)} \left[ R_i^{(v)} - \E\left(R^{(v)}\right) - \frac{\E\left(R^{(u)}R^{(v)}\right) - \E\left(R^{(u)}\right)\E\left(R^{(v)}\right)}{\E\left(R^{(u)}\right)} \right] + o_P(1) \\
    &= \frac{1}{\sqrt{n}} \sum_{i=1}^n \tilde{Y}_i^{(u)} R_i^{(u)} \left( R_i^{(v)} - \frac{\E\left(R^{(u)}R^{(v)}\right)}{\E\left(R^{(u)}\right)} \right) + o_P(1).
\end{align*}

We can now define the asymptotic influence function for any arbitrary pair index as:
\begin{align*}
    \psi_i^{(u,v)} = \left(Y_i^{(u)} - \mu^{(u)}\right) R_i^{(u)} \left( R_i^{(v)} - \frac{\E\left(R^{(u)}R^{(v)}\right)}{\E\left(R^{(u)}\right)} \right).
\end{align*}
By the Multivariate Central Limit Theorem, the joint distribution of these sum statistics converges to a multivariate normal distribution with a covariance matrix defined by $\Omega_{(u,v),(r,s)} = \E\left(\psi_i^{(u,v)} \psi_i^{(r,s)}\right)$ \citep[][Ch. 5 and 7, resp.]{van2000asymptotic, serfling1980approximation}. Given that under MCAR, the data values and the response mechanism metrics are entirely independent, the expectation factors directly into:
\begin{align*}
    \E\left(\psi_i^{(u,v)} \psi_i^{(r,s)}\right) &= \E\left[\left(Y^{(u)} - \mu^{(u)}\right)\left(Y^{(r)} - \mu^{(r)}\right)\right] \\
   & \qquad \qquad \qquad \cdot \E\left[ R^{(u)} R^{(r)} \left( R^{(v)} - \frac{\E(R^{(u)}R^{(v)})}{\E(R^{(u)})} \right) \left( R^{(s)} - \frac{\E(R^{(r)}R^{(s)})}{\E(R^{(r)})} \right) \right],
\end{align*}
which yields exactly equation \eqref{omega_generic}. The final multivariate convergence statement follows directly from applying Slutsky's theorem to the consistent sample-variance plug-in matrix $\hat{\Omega}^{-1}$. This completes the proof.

\subsection{Proof of Lemma \ref{lemma_algorithm}}\label{appendix_lemma_algorithm}

We have that
    \begin{align*}
        \Cov \left( R_1, R_3 \right) &= \E \left( R_1R_3 \right) - \E (R_1)\E (R_3) \\
        &= \E \big( I\left\{ U \leq  {\rho}\right\} R_1 + I\left\{ U > {\rho} \right\}R_1R_2 \big) - \E \left( R_1\right) \E \big( I\left\{ U \leq  {\rho}\right\} R_1 + I\left\{ U > {\rho} \right\}R_2 \big) \\
        &= \Probb \left\{ U \leq {\rho}  \right\} \E (R_1) + \Probb \left\{ U > {\rho} \right\} \E(R_1)\E(R_2) - \E (R_1)\E (R_1) \Probb \left\{ U \leq {\rho}\right\}  - \E (R_1)\E(R_2) \Probb \left\{  U > {\rho} \right\} \\
        &= {\rho}q + (1-{\rho})q^2 -{\rho}q^2  - (1-{\rho})q^2  \\
        &= {\rho}q (1-q).
    \end{align*}
    Additionally,
    \begin{align*}
        \sqrt{\D (R_1)} = \sqrt{q(1-q)}
    \end{align*}
    and
    \begin{align*}
        \D (R_3) &= \E \left(R_3^2\right) - \left(\E (R_3) \right)^2 \\
        &= \E \big( I \left\{ U \leq {\rho}\right\}R_1 + 2 I \left\{ U \leq {\rho} \right\} I \left\{ U > {\rho} \right\} R_1R_2 + I \left\{ U > {\rho} \right\} R_2 \big) - \big( rq + (1-r)q \big)^2 \\
        &=rq +2\E \big( I \left\{ U \leq {\rho}\right\} \left( 1- I \left\{ U \leq {\rho} \right\} \right) \big) q^2 + (1-{\rho})q - q^2 \\
        &=q - q^2 + 2 \cdot 0 \\
        &= q(1-q),
    \end{align*}
    so $\sqrt{\D (R_3)} = \sqrt{q(1-q)}$. Finally, we have that
    \begin{align*}
        \mathbf{Cor} (R_1, R_3) = \frac{\Cov (R_1, R_3)}{ \sqrt{\D (R_1)} \sqrt{\D (R_3)}} = \frac{{\rho}q(1-q)}{\sqrt{q(1-q)}\sqrt{q(1-q)}} = {\rho}.
    \end{align*}

\end{document}